\documentclass{article}
\usepackage{amsmath,amssymb,url,hyperref,pdflscape}
\DeclareMathOperator*{\argmin}{arg\,min}

\usepackage{mathptmx}

\usepackage{graphicx,multirow,subfig} %% ajp removed tikz

\usepackage{color}

\newcommand{\AJPii}[1]{}

\usepackage{mathtools}

\DeclarePairedDelimiter\floor{\lfloor}{\rfloor}

\usepackage{float}
\newfloat{overview}{p}{overview}[section]
\floatname{overview}{Overview}

\usepackage{algorithmic}    % style file for typesetting algorithms
\usepackage{algorithm}      % needed to get the algorithm in a figure
\floatname{algorithm}{Algorithm Schema}

\newenvironment{problem}[1]
    {\begin{center}
    
    \begin{tabular}{|p{0.9\textwidth}|}
    \hline\\
    {#1}\\[1ex]
    }
    { 
    \\\\\hline
    \end{tabular} 
    \end{center}
    }

\newenvironment{sketch}[1][Proof (sketch)]{\begin{trivlist}
\item[\hskip \labelsep {\bfseries #1}]}{\end{trivlist}}

\usepackage{longtable,caption} %ltcaption, caption
\usepackage{apacite} %  

\usepackage[english]{babel}
\usepackage{bbm}

\newcommand{\LB}{\ensuremath {\mathcal X}}
\newcommand{\OUR}{\ensuremath {\mathcal Y}}
\newcommand{\TT}{\ensuremath {\mathcal R}}
\renewcommand{\S}{\ensuremath {\mathcal S}}

\newcommand{\OP}{\ensuremath {\mathcal A}}

\newcommand{\R}{\ensuremath{\mathbbm{R}}}
\newcommand{\N}{\ensuremath{\mathbbm{N}}}

\newcommand{\norm}[1]{\left\| {#1} \right\|}

\newcommand{\bck}[1]{\langle #1\rangle}
\newcommand\scal[2]{\bck{{#1},{#2}}} 

\newtheorem{prop}{Proposition}{\bfseries}{}
{\bfseries}{}

\DeclareMathOperator{\st}{\; s. t. \;}
\DeclareMathOperator{\trace}{trace}
\DeclareMathOperator{\s.t.}{\; s. t. \; }

\DeclareMathOperator{\rank}{rank}
\DeclareMathOperator{\kron}{\,\otimes\,}

\long\def\IGNORE#1{}
\long\def\TODO#1{}
\long\def\NOTE#1{}
\usepackage[T1]{fontenc}
\usepackage[utf8]{inputenc}
\usepackage{authblk}

\title{Semidefinite Programming \\ in Timetabling and Mutual-Exclusion Scheduling}
\author[1]{Jakub~Mare{\v c}ek}
\author[2]{Andrew~J.~Parkes}
\affil[1]{IBM Research -- Ireland}
\affil[2]{The University of Nottingham}

\begin{document}

\maketitle

\begin{abstract}
In scheduling and timetabling applications, the mutual-exclusion constraint stipulates that certain pairs of tasks that cannot be executed at the same time.
This corresponds to the vertex colouring problem in graph theory, for which there are well-known semidefinite programming (SDP) relaxations.
In practice, however, the mutual-exclusion constraint is typically combined with many other constraints, whose SDP representability has not been studied.

We present SDP relaxations for a variety of mutual-exclusion scheduling and timetabling problems, starting from a bound on the number of tasks executed within each period, which corresponds to graph colouring bounded in the number of uses of each colour.
In theory, this provides the strongest known bounds for these problems that are computable to any precision in time polynomial in the dimensions.
In practice, we report encouraging computational results on 
  random graphs, 
  Knesser graphs, ``forbidden intersection'' graphs, the Toronto benchmark, and 
the International Timetabling Competition.
\end{abstract}

%\keywords{ timetabling \and bounded colouring \and graph colouring \and vertex colouring \and semidefinite programming \and augmented Lagrangian \and boundary point method }

\section{Introduction}

Across many areas of combinatorial optimisation, semidefinite programming (SDP) \cite{MR1778223,anjos2011handbook}
 has made it possible to derive strong lower bounds \cite{MR1778234},
  as well as to obtain very good solutions using randomised rounding \cite{lau2011iterative}.
Nevertheless, there seem to be only few applications to practical scheduling, timetabling, or rostering problems.

In scheduling and timetabling problems, one encounters extensions of the  mutual-exclusion constraint, which stipulates that certain pairs of tasks or events  cannot be executed at the same time.
This corresponds to the graph colouring problem in graph theory, for which there are well-known semidefinite programming (SDP) relaxations.
The SDP representability of combinations of the mutual-exclusion constraint with  other constraints has been an open problem.

Clearly, the work on graph colouring provides a test of infeasibility for timetabling and scheduling problems incorporating the mutual exclusion problem. 
Such an infeasibility test compares a lower bound on the optimum of bounded vertex colouring of the conflict graph 
against the number of periods available. 
Lower bounds obtained in ignorance of the extensions, especially the bound on the number of uses of each colour,
  are still perfectly valid, but generally weak. 

In this paper, we set out to explore applicability of semidefinite programming in scheduling and timetabling problems, which extend graph colouring.
In educational timetabling, this corresponds to considering room sizes, room features, room stability, and pre-allocated assignments.
In transportation timetabling, these correspond to considering vehicle capacity, line-vehicle compatibility, etc.
We show that semidefinite programming relaxations of 
 a variety of such problems can be formulated, starting with 
 bounded vertex colouring of the conflict graph.

Our paper is organised as follows: Following some preliminaries, we present our relaxations in Section~\ref{sec:problems}. 
In Section~\ref{sec:algo} we specialise a  well-known first-order method to solving the relaxations and showcase two algorithms for rounding in the relaxations. 
In Section~\ref{sec:analysis}, we analyse some properties of the relaxations and the  performance of algorithms applied to them. 
In Section~\ref{sec:computational}, we present results of extensive computational tests.
On conflict graphs from the International Timetabling Competition 2007, the Toronto benchmark,
as well as on random graphs, we show the relaxations often provide the best possible lower bound and make it possible to 
 obtain very good solutions by randomised rounding.
On ``forbidden intersections'' graphs, we show the strength and weakness of the bound.
Further avenues for research are suggested in Section~\ref{sec:conclusions}.

\section{Notation and Related Work}
\label{sec:related}

%\subsection{The Problems and their Complexity}

\subsection{Semidefinite Programming}
\label{sec:sdp-def}

Firstly, let us revisit the definition of semidefinite programming, 
which is a popular extension of linear programming.
In linear programming (LP), the task is to optimise a linear combination $c^T x$ subject to $m$ linear constraints 
$Ax = b$ subject to the element-wise restriction of variable $x$ to non-negative real numbers. Notice $c$ is an $n$-vector, 
$x$ is a compatible vector variable, $b$ is an $m$-vector, and $A$ is $m \times n$ matrix. 
One can state the problem also as:
\begin{align}
z = \min_{x} c^T x \s.t. \OP_A(x) = b \mbox{ and } x \ge 0 \tag{P LP} \label{lp-p} \\
z = \max_{y} b^T y \s.t. \OP_A^*(y) \le c \tag{D LP} \label{lp-d}  
\end{align}
where linear operator $\OP_A(x)$ (parametrised by matrix $A$) maps vector $x$ to a vector $A x$, and $x \ge 0$ denotes the element-wise non-negativity of $x$, $x \in (\R^+)^n$. 
The element-wise non-negativity of $x$ should be seen as a restriction of vector $x$ to the positive orthant,
which is a convex cone, as all linear combinations with non-negative coefficients of element-wise non-negative vectors 
are element-wise non-negative vectors.
%(Note that the adjoint of $\OP_A(x)$ is $\OP_A^*(y) = A^T y$.)
Using a variety of methods, linear programming can be solved to any fixed precision in polynomial time.
These methods work for other symmetric convex cones as well. 

Semidefinite programming (SDP, \citeNP{Bellman1963,MR1315703,MR1778223,MR1778223,anjos2011handbook}) is a convex optimisation problem,
which generalises linear programming. 
It replaces the vector variable with a square symmetric matrix variable and the polyhedral symmetric convex cone of 
  the positive orthant with the non-polyhedral symmetric convex cone of positive semidefinite matrices. The primal-dual pair in the standard form is:
\begin{align}
z_p =  \min_{X \in \S^n} \scal{C}{X} % + f^T t
\s.t.\; \OP_A(X) = b \mbox{ and } X \succeq 0 \tag{P STD} \label{std-p} \\  
z_d = \max_{y \in \R^m, S \in \S^n} b^T y 
\s.t. \; \OP_A^*(y) + S = C \mbox{ and } S \succeq 0 \label{std-d} \tag{D STD}  
\end{align}
where $X$ is a primal variable in the set of $n \times n$ symmetric matrices $\S^n$, $y$ and $S$ are the corresponding dual variables,
$b$ is an $m$-vector, $C$, $A_i$ are compatible matrices, and linear operator $\OP_A(X)$ (parametrised by a symmetric matrix $A$) maps a symmetric matrix $X$ to vectors in $\R^m$,
wherein the $i$th element $\OP_A(X)_i = \scal{A_i}{X}$. % = \sum_{i,j} A_{i,j} X_{j,i}$,
%linear operator $\OP_B(X)$ maps symmetric matrices to vectors $\R^q$, where $i$th element $B(X)_i = \scal{B_i}{X}$, % = \sum_{i,j} B_{i,j} X_{j,i}$,
$\OP_A^*(y)$ is again the adjoint operator.
%$M \ge 0$ denotes $M$ is element-wise non-negative,
$M \succeq N$ or $M - N \succeq 0$ denotes $M - N$ is positive semidefinite.
Note that an $n \times n$ matrix, $M$, is positive semidefinite if and only if $y^T M y \ge 0$ for all $y \in \R^n$. 
One can also treat inequalities explicitly in the primal-dual pair: 
\begin{align}
%\min_{X,t} \scal{C}{X} + f^T t \s.t. X \succeq 0 \mbox{ and } \scal{A_i}{X} + G_i^T t = b_i \tag{P} \label{sdp-p} \\
%\max_{y} b^T y \s.t. C - \sum_{i} A_i y_i \succeq 0 \mbox{ and } \sum{i} G_i y_i = f \label{sdp-d} \tag{D}
z_p =  \min_{X \in \S^n} \scal{C}{X} % + f^T t
\s.t.\; \OP_A(X) = b \mbox{ and } \OP_B(X) \ge d \mbox{ and } X \succeq 0 
\tag{P SDP} \label{sdp-p} \\  
z_d = \max_{y \in \R^m, v \in \R^q, S \in \S^n} b^T y + d^T v
\s.t. \; \OP_A^*(y) + \OP_B^*(v) + S = C \mbox{ and } S \succeq 0 \mbox{ and } v \ge 0 
\label{sdp-d} \tag{D SDP}  
\end{align}
where $d$ is a $q$-vector and linear operator $\OP_B(X)$ maps $n \times n$ matrices to $q$-vectors similarly to $\OP_A$ above.
As all linear combinations with non-negative coefficients of positive semidefinite matrices 
are positive semidefinite, $X \succeq 0$ should again be seen as a restriction to a convex cone.
%Unlike in linear programming, strong duality $z_p = z_d$ holds only if certain conditions hold; One such condition, Slater's, requires both primal (\ref{sdp-p}) and dual (\ref{sdp-d}) problems to be feasible. Otherwise, the similarity of the primal-dual pairs in linear (\ref{lp-p}, \ref{lp-d}) and semidefinite (\ref{sdp-p}, \ref{sdp-d}) programming are clear.

\subsection{Semidefinite Programming Relaxations of Graph Colouring}

Graph colouring, also known as vertex colouring, or partition into independent sets, is the problem of: 

\begin{problem}{\sc Graph Colouring}
  Given an undirected graph $G = (V, E)$ with vertices $V = {1, 2, \ldots, n}$ and edges $E \subset \{ (u, v) \s.t. 1 \le u < v \le n \}$,
  return a partition $P = (P_i)$ of $V$ of the smallest possible cardinality $|P|$ so that for each partition $P_i$, for no edge $(u, v) \in E$, there are both $u$ and $v$ in $P_i$.
  As in any partition, $\cup_i P_i = V$ and for all $1 \le i < j \le |P|$, we have $P_i \cap P_j = \emptyset$.
\end{problem}

The partitions are known as ``colour classes'' or ``independent sets'' in graph colouring, or (assignment to) ``time periods'' in timetabling and scheduling.
The optimum, i.e., the smallest possible number $|P|$ of colour classes, is denoted the ``chromatic number'' in graph colouring
or minimum number of required time periods in timetabling, 
or ``makespan'' in connection with certain mutual-exclusion problems (cf. Section \ref{sec:prob:mutualex}) in  
the scheduling literature. 

The decision version of graph colouring appears on Karp's original list \cite{MR0378476} of NP-Complete problems.
In polynomial time, one can obtain lower bounds on the chromatic number,
 for instance using linear or semidefinite programming.
Just as there are a number of ways of formulating a lower bound on the chromatic number in linear programming,
there are a number of ways of formulating a lower bound on the chromatic number using SDP.
All are, in some sense, related to the inequality of Wilf \cite{wilf1967eigenvalues},
wherein the largest eigenvalue of an adjacency matrix of a graph incremented by one bounds the chromatic number of a graph from the above. 
By considering the semi-definite programming lower bound on the largest-eigenvalue upper-bound, one obtains a parameter of a graph, sometimes known as `theta', $\theta$, \cite{Lovasz1979}, which is at least as large as the clique number and no more than the chromatic number, yet is computable in polynomial time using SDP.
The known bounds for colouring form a hierarchy \cite{1399025}: 
\begin{align}
\alpha(G)   \le \LB'(G)             \le \LB(G)             \le \LB^+(G)             \le \LB^{+\bigtriangleup}(G)             \le \chi(\overline G), & \mbox{ or } \hfill \\  
\omega(G)   \le \LB'(\overline G)   \le \LB(\overline G)   \le \LB^+(\overline G)   \le \LB^{+\bigtriangleup}(\overline G)   \le \chi(G), & \notag
\end{align}
where $\alpha$ is the size of the largest independent set, 
$\omega$ is the size of the largest clique,
$\chi$ is the chromatic number,
$\LB(G) = \theta(\bar G)$ is the vector chromatic number \cite{Lovasz1979,Karger1998}, 
$\LB'(G) = \theta_{1/2}(\bar G)$ is the strict vector chromatic number \cite{Kleinberg1998}, 
$\LB^+(G) = \theta_{2}(\bar G)$ is the strong vector chromatic number \cite{1399025}, and bar indicates complementation of a graph.
%Also (Coja-Oghlan, 2003 ***), for $c_0/n \le p \le 1/2$:
%$$\prob(c_1 \sqrt{np} \le \theta_{1/2}(G_{n,p}) \ge 1 - exp(-n)$$
% $$\prob(\theta_{1/2}(G_{n,p} \le c_2 \sqrt{np}) \ge 1 - exp(-n)$$ ***
%where $c0, c_1$, and $c_2$ are constants.
%Typically, the analysis of the lower bounds use the vector programming formulation,
%which is possible to reformulate as primal (\ref{sdp-p}) and dual (\ref{sdp-d}) semidefinite programming problems.
In Figure~\ref{over1}, we summarise all three formulations for all three lower bounds in the vector programming notation.
None of the formulatons has, however, been extended to bounded colouring, up to the best of our knowledge.

\begin{landscape}
\begin{figure}
\caption{An overview of known vector programming (VP) and semidefinite programming (SDP) relaxations of vertex colouring,\\ 
$\alpha(G)\le \LB'(G) \le \LB(G) \le \LB^+(G) \le \chi(\overline G)$ or 
$\omega(G)\le \LB'(\overline G) \le \LB(\overline G) \le \LB^+(\overline G) \le \chi(G)$. }
\label{over1}

\begin{minipage}[t]{0.3\linewidth}\centering
Lov{\'a}sz's Bound as VP\\
  \begin{align}
 \LB(G) = & \min t = \theta(\bar G) \label{theta-vp1} \\
\st: & \norm{v_i}_2 = 1 & \forall i \in V \notag \\
     & \scal{v_i}{v_j} = - \frac{1}{t - 1} \quad \forall \{i,j\} \in E \notag
  \end{align}
\end{minipage}
\hspace{0.5cm}
\begin{minipage}[t]{0.3\linewidth}
\centering
Primal SDP for Lov{\'a}sz's Bound\\
  \begin{align}
\LB(G) = & \max \scal{J}{X} = \theta(\bar G)\label{theta-sdp1p} \\
\st: & \trace(X) = 1   \notag \\
     & X_{uv} = 0 \quad \forall \{u, v\} \in E \notag \\
     & X \succeq 0 \notag
  \end{align}
\end{minipage}
\hspace{0.5cm}
\begin{minipage}[t]{0.3\linewidth}
\centering
Dual SDP of Lov{\'a}sz's Bound\\
  \begin{align}
\LB(G) = & \min t = \theta(\bar G)\label{theta-sdp1d} \\
\st: & U_{uu} = 1 \quad \forall u \in V  \notag \\
     & U_{uv} = -\frac{1}{t-1} \quad \forall \{u, v\} \in \bar E \notag \\
     & U \succeq 0, t \ge 2 \notag
  \end{align}
\end{minipage}

\vskip 6mm

\begin{minipage}[t]{0.3\linewidth}\centering
Kleinberg's Bound as VP\\
  \begin{align}
\LB'(G) = & \min t = \theta_{1,2}(\bar G) \label{theta-vp2} \\
\st: & \norm{v_i}_2 = 1 & \forall i \in V \notag \\
     & \scal{v_i}{v_j} \le - \frac{1}{t - 1} \forall \{i,j\} \in E \notag
  \end{align}
\end{minipage}
\hspace{0.5cm}
\begin{minipage}[t]{0.3\linewidth}
\centering
Primal SDP for Kleinberg's Bound\\
  \begin{align}
\LB'(G) = & \max \scal{J}{X} = \theta_{1,2}(\bar G)\label{theta-sdp2p} \\
\st: & \trace(X) = 1   \notag \\
     & X_{uv} = 0 \quad \forall \{u, v\} \in E \notag \\
     & X_{uv} \ge 0 \quad \forall \{u, v\} \in \bar E \notag \\
     & X \succeq 0 \notag
  \end{align}
\end{minipage}
\hspace{0.5cm}
\begin{minipage}[t]{0.3\linewidth}
\centering
Dual SDP of Kleinberg's Bound\\
  \begin{align}
\LB'(G) = & \min t = \theta_{1,2}(\bar G)\label{theta-sdp2d} \\
\st: & U_{uu} = t \quad \forall u \in V \notag \\
     & U_{uv} \le - \frac{1}{t - 1} \quad \forall \{u, v\} \in \bar E \notag \\
     & U \succeq 0, t \ge 2 \notag
  \end{align}
\end{minipage}

\vskip 6mm

\begin{minipage}[t]{0.3\linewidth}\centering
Szegedy's Bound as VP
  \begin{align}
\LB^+(G) = & \min t = \theta_2(\bar G)\label{theta-vp2-2} \\
\st: & \norm{v_i}_2 = 1 & \forall i \in V \notag \\
     & \scal{v_i}{v_j} = - \frac{1}{t - 1} \quad \forall \{i,j\} \in E \notag \\
     & \scal{v_i}{v_j} \geq - \frac{1}{t - 1} \quad \forall \{i,j\} \in \bar E \notag 
  \end{align}  
\end{minipage}
\hspace{0.5cm}
\begin{minipage}[t]{0.3\linewidth}
\centering
Primal SDP of Szegedy's Bound\\
  \begin{align}
\LB^+(G) = & \max \scal{J}{X} = \theta_2(\bar G) \label{theta-sdp3p} \\
\st: & \trace(X) = 1   \notag \\
     & X_{uv} \le 0 \quad \forall \{u, v\} \in E \notag \\
     & X \succeq 0 \notag
  \end{align}    
\end{minipage}
\hspace{0.5cm}
\begin{minipage}[t]{0.3\linewidth}
\centering
Dual SDP of Szegedy's Bound\\
  \begin{align}
         \LB^+(G) = & \min t = \theta_2(\bar G)\label{theta-sdp3d} \\
\st: & U_{uu} = t \quad \forall u \in V \notag \\
     & U_{uv} = - \frac{1}{t - 1} \quad \forall \{u, v\} \in \bar E \notag \\
     & U_{uv} \ge - \frac{1}{t - 1} \quad \forall \{u, v\} \in E \notag \\
     & U \succeq 0, t \ge 2 \notag
  \end{align}
\end{minipage}

\end{figure}
\end{landscape}

There are a number of ways of deriving and thinking about the SDP relaxations.
In any case, the primal $n \times n$ matrix variable can be seen  
% in the semidefinite programming formulations (\ref{theta-sdp1p,theta-sdp2p,theta-sdp3p}) 
as
\begin{align}
M_{u,v} = \begin{cases}
\; 1 & \text{ if vertex $u$ is in the same colour class as $v$ } \\
\; 0 & \text{ otherwise. }
\end{cases}
\end{align}

Notice matrix $M$ has the ``hidden block diagonal''property:

\begin{prop}
For any value of $M$, there exists a permutation matrix $P$, such that
\begin{align}
M^{bd} = P^T MP = 
\left[ {\begin{aligned}
   J_{c_1} \; & \; 0 \; & {} & 0  \\ 
   0 \; & \; J_{c_2} \; & {} & 0  \\
   {} & {} &  \ddots  & {}  \\ 
   0 \; & \; 0 \; & {} & \; J_{c_s}  
 \end{aligned} } \right]
\label{eqn:blockdiagonal}
\end{align}
where $J_c$ is the $c \times c$ matrix of all ones and $\sum_{i =1}^{s} c_i = n$. Such $P^T MP$ is denoted a direct sum of $J_c$. 
\end{prop}

One can hence derive the semidefinite programming relaxation from:

\begin{prop}[\citeA{Dukanovic2004-patat}]
For any symmetric 0-1 matrix $M$ there exists a permutation matrix $P$ such that $P^T M P$ is the direct sum of $s$ all-ones matrices
if and only if there is a vector of all-ones on the diagonal and the rank of $M$ is $s$ and $M$ is positive semidefinite.
\end{prop}

by relaxing the rank constraint, as usual \cite{Fazel2004}.
Alternatively, one can see theta as an eigenvalue bound, where the largest eigenvalue $\lambda_{\max}(A) = \min \{ t \st tI - A \succeq 0 \}$ for an identity matrix $I$.
Perhaps most ``fundamentally'', one could see theta as a relaxation of the co-positive programming formulation 
of graph colouring, recently proposed by \citeA{Bomze2010}.
A number of other derivations have been surveyed by \citeA{Knuth1994}. 

\subsection{Related Applications}

Within the job-shop scheduling, SDP relaxations of the maximum cut problem (MAXCUT, \citeNP{MR1412228}) have been adapted to 
  scheduling workload on two machines \cite{MR1868715,MR1999701} and home-away patterns in sports scheduling \cite{MR2294048}.
%Still restricted to two machines, \cite{MR1999701} introduced additional constraints and improved the approximation ratio for a procedure rounding the relaxation. 
\cite{Bansal2016} extended this to scheduling with weighted completion time objectives on any number of unrelated machines, 
 using a clever rounding technique in a lift-and-project relaxations.
They have also shown that the relaxation of \cite{MR1868715} 
is in some sense weak (has 3/2 integrality gap).
We are not aware of any applications of semidefinite programming to mutual-exclusion scheduling or timetabling, excepting two abstracts of the present authors \cite{Marecek2010-patat,Marecek2012-patat}, which we build upon in this paper.

\section{Old Problems and Novel Relaxations}
\label{sec:problems}

\subsection{Mutual-Exclusion Scheduling}
\label{sec:prob:mutualex}

In one of the prototypical problems in timetabling, scheduling, and staff rostering \cite{Welsh1967,Burke2004,MR1403900}, 
one needs to assign $n$ unit-time events, classes, tasks, or jobs (``vertices''),
some of which must not be run, executed, or taught at the same time (``the mutual-exclusion constraint''), perhaps due to the use to some shared, renewable resource, 
to $m$ rooms, processors, machines, or employees (``uses of a colour''),
so that the number of units of time required (``makespan'', ``number of colours'') is as small as possible. 
It is natural to represent the elements being assigned by the elements of set $V = {1, 2, \ldots, n}$,
and to represent the mutual-exclusion constraint by a set of pairs of elements of $V$, 
which are not to be assigned to disjoint time-units (``conflict graph'').

Let us now consider:

\begin{problem}{\sc $m$-Bounded Colouring}
  Given an undirected graph $G = (V, E)$ with vertices $V = {1, 2, \ldots, n}$ and edges $E \subset \{ (u, v) \s.t. 1 \le u < v \le n \}$,
  and an integer $m \le |V|$,
  return a partition $P = (P_i)$ of $V$ of the smallest possible cardinality $|P|$ so that for each partition $P_i$, $|P_i| \le m$ and for no edge $(u, v) \in E$ there are both $u$ and $v$ in $P_i$.
\end{problem} 
  
Terminology varies. 
In the scheduling community, the problem is known as 
scheduling of unit-length tasks on $m$ parallel machines with renewable resources \cite[Chapter 12]{blazewicz2007handbook}, 
{\sc Mutual-Exclusion Scheduling} \cite{MR1403900}, 
scheduling with incompatible objects \cite{leighton1979graph}, or P$m$ | res,$p_j=1$ | C$_{\max}$ in the notation of \cite{graham1979optimization}. 
In discrete mathematics \cite{KRARUP19821}, some authors \cite{MR1265071} refer to the problem as {\sc Partition into Bounded Independent Sets}, while others use {\sc $m$-Bounded Colouring} \cite{MR1210106} or just Bounded Graph Colouring.
There is a simple transformation to {\sc Equitable Colouring} \cite{MR2183465}, 
 where the cardinality of colour classes can differ by at most one.
In terms of complexity theory, $m$-Bounded Colouring is in P on trees \cite{MR1275266}, 
but NP-Hard on co-graphs, interval graphs, and bipartite graphs \cite{MR1265071}.

From the block-diagonal property of the matrix variable (\ref{eqn:blockdiagonal}), it seems clear that we expect row-sums and column-sums in the binary-valued variable to be bounded from above by $m$:

  \begin{align}
             & \min \; t  \\
\st: \qquad\qquad
& X_{vv} = 1 \quad \forall v \in V  \label{eqn:X-linear-eq1} \tag{E1} \\
     & X_{uv} = 0 \quad \forall \{u, v\} \in E \label{eqn:X-linear-eq2} \tag{E2} \\
     & \rank(X) = t \label{eqn:rankt} \tag{rank-$t$} \\
     & \sum_{u \in V } X_{uv} \le m \quad \forall v \in V  \label{eqn:linear-ineq1} \tag{IN} \\
& X\succeq 0 \label{eqn:X-psd} \tag{PSD} \\
     & X_{uv} \in \{ 0, 1 \} \label{eqn:X-range} \tag{Binary}
\end{align}

where $X$ is an $n \times n$ matrix variable and $t$ is a scalar variable.

By dropping the element-wise integrality, 
relaxing a non-convex bound on rank($X$) to a convex bound on the trace($X$), and algebraic manipulations,
we obtain the relaxation:

  \begin{align}
             & \min \; t&  \label{eqn:OUR-SDP-Simple} \\
\st \qquad Y_{vv} & = t & \quad \forall v \in V  \label{eqn:linear-eq1}  \tag{E1} \\
     Y_{uv} & = 0 & \quad \forall \{u, v\} \in E \label{eqn:linear-eq2} \tag{E2} \\
     \sum_{u \in V } Y_{uv} & \le tm & \quad \forall v \in V  \tag{IN} \\
     Y - J & \succeq 0 \label{eqn:psd} \tag{PSD} & \\
     Y_{uv} & \ge 0 & \quad \forall \{u, v\} \label{eqn:nonneg-bound} \tag{N} %\\
\end{align}
where $Y$ is an $n \times n$ matrix variable and $t$ is a scalar variable.
This can be seen as a spectahedron (\ref{eqn:linear-eq1}--\ref{eqn:psd}) being intersected by a polyhedron given by 
the linear inequalities (\ref{eqn:linear-ineq1}).

Notice that additional inequalities $\sum_{v \in V } Y_{uv} \le tm \quad \forall u \in V$ are not required due to symmetry.
Notice also that the usual theta-like relaxations ($\LB(\overline G), \LB'(\overline G), \LB^+(\overline G)$) cannot be used easily, 
as one cannot easily work with a graph's complement, both due to its density and due to the constraints on the representability of our extensions.

The complication is the relaxation above is not a semidefinite program in the standard form (\ref{std-p}, \ref{std-d}).
Notice scalar variable $t$ has been introduced only for clarity.
As long as the entries on the diagonal of the matrix variable are constrained to be equal, any one of them can be used instead.
One can either introduce new scalar slack variables and convert inequalities to equalities,
 or one can design solvers treating inequalities explicitly.

There remains the constraint $Y - J \succeq 0$ to deal with, 
as the standard form only allows to require a matrix, rather than expression, to be psd.
The mechanistic approach, employed by automated model transformation tools \cite{YALMIP}, for instance,
is to double the dimension by introducing new variable $X$, set $Y - J = Z, Z \succeq 0$. For an arbitrary vertex $w \in V$, one obtains: 
  \begin{align}
     & \max Y_{w,w} \label{eqn:Rewritten} \\
\st \qquad & Z_{u,v} = -1 \quad \forall \{u, v\} \in E \label{eqn:A1} \tag{A1} \\
     & Z_{v,v} = Z_{w,w} \quad \forall v \in V \setminus \{ w \} \label{eqn:A2} \tag{A2} \\
     & Y_{u,v} - Z_{u,v} = 1 \quad \forall u, v \in V  \label{eqn:W} \tag{W} \\
     & \sum_{u \in V } Y_{uv} \le tm \quad \forall v \in V  \label{eqn:B} \tag{B} \\
     & Z \succeq 0  \notag
\end{align}

An alternative approach is to optimise:
  \begin{align}
     & \min X_{w,w} + 1\label{eqn:Scaled} \\
\st \qquad & X_{u,v} = -1 \quad \forall \{u, v\} \in E \tag{A1} \\
     & X_{v,v} = X_{w,w} \quad \forall v \in V \setminus \{ w \} \tag{A2} \\
     & |V| - m X_{w,w} - m - 1 + \sum_{u \in V } X_{uv} \le 0 \quad \forall v \in V \tag{B} \\
     & X \succeq 0  \notag
\end{align} 
Any formulation involving $Y - J \succeq$ in this paper can be easily transformed in this fashion.

\subsection{Initial Assignment}

In many applications, one has to deal with complicating constraints. 
In graph-theoretic terms, 
the most common complicating constraints are pre-existing assignments.  
Pre-existing assignments can be represented as subsets of $V$, which need to be assigned to the same unit of time.    
This corresponds to:

\begin{problem}{\sc $m$-Bounded Colouring with Pre-Colouring}
  Given an undirected graph $G = (V, E)$ with vertices $V = {1, 2, \ldots, n}$ and edges $E \subset \{ (u, v) \s.t. 1 \le u < v \le n \}$,
  an integer $m \le |V|$,
  and a family $C = (C_i)$ of disjoint subsets of $V$ with $|C_i| \le m$, 
  return a partition $P = (P_i)$ of $V$ of the smallest possible cardinality $|P|$ so that for each partition $P_i$, $|P_i| \le m$ and for no edge $(u, v) \in E$ there are both $u$ and $v$ in $P_i$
  and for each set $C_i \in C$, there is a partition in $C_i \subseteq P_j \in P$.
  The graph $G$ is called ``conflict graph'' and family $C$ is called ``pre-colouring''.
  The partition $P$ corresponds to the ``same-colour equivalence'' and $P_i$ are called ``colour classes'' or ``independent sets''.
\end{problem}

In terms of complexity theory, $m$-Bounded Colouring with Pre-Colouring is NP-Hard even on trees \cite{MR2183465}.
Even in the case of trees, however, there are fixed parameter tractable algorithms \cite{MR1265071,MR2183465}.

In terms of SDP representability, 
given a pre-colouring $C = (C_i), C_i \subseteq V$, it suffices to set the corresponding elements of matrix variable $Y$ to $t$:
  \begin{align}
\OUR(G, m) = & \max t \label{OUR-SDP} \\
\st: & Y_{vv} = t \quad \forall v \in V  \label{eqn:linear-eq1-2} \tag{E1}  \\
     & Y_{uv} = 0 \quad \forall \{u, v\} \in E \tag{E2} \\
     & Y_{uv} = t \quad \forall u, v \in C_i, u \ne v, C_i \in X \label{eqn:linear-eq3-2}  \tag{E3} \\
     & Y_{uv} = 0 \quad \forall u \in C_i \in C, v \in C_j \in C, C_i \ne C_j \label{eqn:linear-eq4-2}  \tag{E4} \\    
     & \sum_{u} Y_{uv} \le tm \quad \forall v \in V  \label{eqn:ineq1} \tag{L1}  \\
     & \sum_{v} Y_{uv} \le tm \quad \forall u \in V  \label{eqn:ineq2} \tag{L2}  \\
     & Y_{uv} \ge 0 \quad \forall u, v \in V, u \ne v, \{u, v\} \not\in E \tag{L3} \\
     & \scal{J}{X} - n m t \le 0 \label{eqn:weird} \tag{L4} \\
     & Y - J \succeq 0. \notag
\end{align}

\subsection{A Reformulation}

As an aside, notice that {\sc $m$-Bounded Colouring with Pre-Colouring} can be easily transformed into a problem without pre-colouring, but with certain 
weights on vertices:   

\begin{problem}{\sc $c$-Weighted $m$-Bounded Colouring}
  Given an undirected graph $G = (V, E)$ with vertices $V = {1, 2, \ldots, n}$ and edges $E \subset \{ (u, v) \s.t. 1 \le u < v \le n \}$,
  a vector of positive integers $c$ of dimension $|V|$,  
  and an integer $m \le |V|$, 
  return a partition $P = (P_i)$ of $V$ of the smallest possible cardinality $|P|$ so that for each partition $P_i$, $c p_i \le m$, where $p_i$ is the $0-1$ index vector corresponding to $P_i$, and for no edge $(u, v) \in E$ there are both $u$ and $v$ in $P_i$.
\end{problem}

For any non-empty $C$, this leads to a reduction in the dimension of the matrix variable, 
compared to the simple relaxation (\ref{OUR-SDP}):

  \begin{align}
\dot \OUR(G, m) = & \max t \label{OUR-SDP-Weighted} \\
\st \qquad & Y_{vv} = t \quad \forall v \in V  \tag{E1}  \\
     & Y_{uv} = 0 \quad \forall \{u, v\} \in E \tag{E2} \\
     & \sum_{u \in V} c_u Y_{uv} \le t m \quad \forall v \in V  \tag{L1}  \\
     & \sum_{v \in V} c_v Y_{uv} \le t m \quad \forall u \in V  \tag{L2}  \\
%     & Y_{uv} \ge 0 \quad \forall u, v \in V, u \ne v, \{u, v\} \not\in E \tag{L3} \\
%     & \scal{J}{X} - nmt \le 0 \\
     & Y - J \succeq 0 \notag
\end{align}

While the reduction improves computational performance, when used with off-the-shelf solvers, it may render both the design of custom solvers and their analysis (cf. Section \ref{sec:analysis}) more challenging.

\subsection{Simple Laminar Timetabling}

In timetabling applications, it is often necessary to consider assignment of events to both periods and rooms.
In the relaxations above, vertices within a single colour class correspond to events taking place at the same time,
but rooms are represented only by the $m$-bound, which corresponds to the number of room available.
This may be insufficient: Consider, for example, a situation with two large lecture rooms, twenty periods per week, 
and forty large lectures.
One could formulate this problem as a binary linear program with a variable with three indices (events, rooms, and periods)
and apply the operators of \citeA{MR1098425} or \citeA{MR2041934} to obtain semidefinite programming relaxations.
We present a number of alternatives, where the matrix variable has a considerably lower dimension.

In particular, one can extend {\sc $m$-Bounded Colouring with Pre-Colouring} to consider not only the number $m$ of available rooms (processors, machines, employees, or similar), 
but also their capacities and features.
This corresponds to: 

\begin{problem}{\sc Simple Timetabling}
  Given an undirected conflict graph $G = (V, E)$ where vertices $V$ are also denoted events,
  a vector of capacity-requirements $p \in \R^{|V|}$,  
  an integer $m \le |V|$,
  a vector of capacities $r \in \R^{m}$,
  a number of features $f \in \N$,  
  set $F \subseteq V \times F$ detailing feature-requirements of events,
  set $G \subseteq \{1, 2, \ldots, m \} \times F$ detailing feature availability,
  and a family $C = (C_i)$ of disjoint subsets of $V$ with $|C_i| \le m$,  
  return a partition $P = (P_i)$ of $V$ of the smallest possible cardinality $|P|$ so that
  \begin{itemize} 
  \item for each partition $P_i$, $|P_i| \le m$ 
  \item for each partition $P_i$ and for no edge $(u, v) \in E$ there are both $u$ and $v$ in $P_i$
  \item for each partition $P_i$ and for each distinct capacity $c$ in $r$,
        the subset of $P_i$ with capacity greater or equal than $c$, according to $p$,
        is less than the number of elements in $r$ greater or equal to $c$
  \item for each partition $P_i$ and for each feature $1, 2, \ldots, f$,
        the subset of $P_i$ requiring the feature, according to $F$, 
        is less than the number of rooms where it is available, according to $G$
  \item for each set $C_i \in C$, there is a partition in $C_i \subseteq P_j \in P$.
  \end{itemize}
\end{problem}

In an important special case, which we denote {\sc Simple Laminar Timetabling}, 
sets of events and rooms with certain feature-requirements and feature-availability are ``laminar''. 
A collection of sets $F$ is called ``laminar'' if $A, B \in F$ implies that $A \subseteq B, B \subseteq A$ or $A \cap B = \emptyset$.
The subsets of rooms requiring capacities larger than a certain value are naturally laminar, 
but naturally occurring features need not be.
As $m$-Bounded Colouring with Pre-Colouring is a special case of Simple Laminar Timetabling,
hardness results cited above apply also to Simple (Laminar) Timetabling.
Such laminar timetabling and associated bounds can be of  interest, for example, in planning the capacities of rooms  cf. \cite{BeyrouthyEtal2010:space-profiles}.

Initially, we restrict ourselves to the Simple Laminar Timetabling and extend the $c$-weighted relaxation above (\ref{OUR-SDP-Weighted}). 
Let us suppose $n$ vertices $V$ of the conflict graph correspond to $n$ events attended
by $p_1, p_2, \ldots, p_n$ persons each, whereas the $m$-bound corresponds to rooms
of capacities $r_1, r_2, \ldots, r_m$. Let us denote distinct 
distinct numbers of persons attending an event $P$,
 $L(p)$ the events of size $s > p$ and
 $R(p)$ the rooms of capacity $r > p$.
Clearly, one can add two constraints (\ref{eqn:LPR1}, \ref{eqn:LPR2}) for each element of $P$:

  \begin{align}
\dot \TT(G, p, r) = & \max t \label{OUR-SDP-Timetabling} \\
\st \qquad & Y_{vv} = t \quad \forall v \in V \notag  \\
     & Y_{uv} = 0 \quad \forall \{u, v\} \in E \notag \\
     & \sum_{u \in L(p) } Y_{uv} \le t |R(p)| \quad \forall p \in P \quad \forall v \in L(p) \label{eqn:LPR1} \tag{PR1} \\
     & \sum_{v \in L(p) } Y_{uv} \le t |R(p)| \quad \forall p \in P \quad \forall u \in L(p) \label{eqn:LPR2} \tag{PR2} \\
     & Y - J \succeq 0 \notag
\end{align}  

Notice the newly added constraints (\ref{eqn:LPR1}, \ref{eqn:LPR2}) subsume the linear inequalities
(\ref{eqn:ineq1}, \ref{eqn:ineq2}) of the relaxations above.
There is always a feasible solution,
and as will be shown in Section~\ref{sec:recovery}, one can efficiently produce 
feasible timetables satisfying those constraints from the value of the matrix variable.

Considering there is always a permutation matrix such that the matrix variable is block-diagonal (\ref{eqn:blockdiagonal}),
one can also add bounds obtained by counting arguments. The simple counting bound is $\sum_{u, v \in V}  Y_{uv} \le m |V|$.
Indeed, there are at most $|V| / m$ blocks of $m^2$ non-zeros each. 
One can generalise the bound to subsets $P$ of events:
  \begin{align}
\ddot \TT(G, p, r) = & \max t \label{OUR-SDP-WithCounting} \\
\st \qquad & Y_{vv} = t \quad \forall v \in V \notag  \\
     & Y_{uv} = 0 \quad \forall \{u, v\} \in E \notag \\
     & \sum_{u \in L(p)} Y_{uv} \le t |R(p)| \quad \forall p \in P \quad \forall v \in L(p) \label{eqn:count:PR1} \tag{PR1} \\
     & \sum_{v \in L(p)} Y_{uv} \le t |R(p)| \quad \forall p \in P \quad \forall u \in L(p) \label{eqn:count:PR2} \tag{PR2} \\
     & \sum_{u \in L(p)} \sum_{ v \in V } Y_{uv} \le m t |R(p)| \quad \forall p \in P  \label{eqn:count:CB1} \tag{CB1} \\
     & \sum_{v \in L(p)} \sum_{ u \in V } Y_{uv} \le m t |R(p)| \quad \forall p \in P  \label{eqn:count:CB2} \tag{CB2} \\     
     & Y - J \succeq 0 \notag
\end{align}  

In {\sc Simple Laminar Timetabling}, one can include feature considerations similarly to capacity considerations. In a slight abuse of notation,
we use $F(f)$ to denote the set of events requiring feature $f$ and $G(f)$ the set of rooms with feature $f$. 
  \begin{align}
\TT(G, p, r, f_{\max}, F, G) = & \max t \label{OUR-SDP-WithFeatures} \\
\st \qquad & Y_{vv} = t \quad \forall v \in V \notag  \\
     & Y_{uv} = 0 \quad \forall \{u, v\} \in E \notag \\
     & \sum_{u \in L(p) } Y_{uv} \le t |R(p)| \quad \forall p \in P \quad \forall v \in L(p) \label{eqn:feat:PR1} \tag{PR1} \\
     & \sum_{v \in L(p) } Y_{uv} \le t |R(p)| \quad \forall p \in P \quad \forall u \in L(p) \label{eqn:feat:PR2} \tag{PR2} \\
     & \sum_{u \in L(p)} \sum_{ v \in V } Y_{uv} \le m t |R(p)| \quad \forall p \in P  \label{eqn:feat:CB1} \tag{CB1} \\
     & \sum_{v \in L(p)} \sum_{ u \in V } Y_{uv} \le m t |R(p)| \quad \forall p \in P  \label{eqn:feat:CB2} \tag{CB2} \\     
     & \sum_{u \in F(f)} Y_{uv} \le t |G(f)| \quad \forall 1 \le f \le f_{\max} \quad \forall v \in F(f) \label{eqn:feat:FR1} \tag{FR1} \\
     & \sum_{v \in F(f)} Y_{uv} \le t |G(f)| \quad \forall 1 \le f \le f_{\max} \quad \forall u \in F(f) \label{eqn:feat:FR2} \tag{FR2} \\
     & \sum_{u \in F(f)} \sum_{ v \in V} Y_{uv} \le m t |G(f)| \quad \forall 1 \le f \le f_{\max} \label{eqn:feat:FC1} \tag{FC1} \\
     & \sum_{v \in F(f)} \sum_{ u \in V} Y_{uv} \le m t |G(f)| \quad \forall 1 \le f \le f_{\max} \label{eqn:feat:FC2} \tag{FC2} \\
     & Y - J \succeq 0 \notag
\end{align}
Notice, however, the limitation to the laminar special case. 

\subsection{Simple Timetabling}

A natural approach to formulating the Simple Timetabling Problem without laminarity requirements uses the additional variables:

\begin{align}
Z_{u,v,r} = \begin{cases}
\; 1 & \text{ if vertex $u$ is in the same colour class as $v$ and $v$ is assigned (room) $r$} \\
\; 0 & \text{ otherwise. }
\end{cases}
\end{align}  

The additional constraints follow. In timetabling language: an event $v$ is assigned a room in $Z_{u,v,r}$ for some $r$
 if and only if it is assigned time in $Y_{u,v}$ (\ref{eqn:rooms1:iff}),
 no event is assigned to two rooms (\ref{eqn:rooms1:exclusive}),
 no two events share a room (\ref{eqn:rooms1:exclusive2})
 and $Z_{u,v,r} = 0$ if the event-room combination does not match the event's room-feature or capacity requirements (\ref{eqn:rooms1:capacity},\ref{eqn:rooms1:feature}).
\begin{align}
\sum_{1 \le r \le m} Z_{u,v,r} = Y_{u,v}     & \quad \forall u, v \in V \label{eqn:rooms1:iff} \\ 
Z_{u,v,r} + Z_{u,v,r'} \le t                     & \quad \forall u,v \in V \forall 1 \le r \le m \forall r' \neq r \label{eqn:rooms1:exclusive} \\
Z_{u,v,r} + Z_{u,w,r} \le t                      & \quad \forall u,v \in V \forall w \in V \setminus \{v\} \label{eqn:rooms1:exclusive2} \\
Z_{u,v,r} = 0                                    & \quad \forall v \in V \quad \forall 1 \le r \le m, p_v \ge r_r \label{eqn:rooms1:capacity}  \\
Z_{u,v,r} = 0                                    & \quad \forall v \in V  \quad \forall 1 \le f \le f_{\max}, (v, f) \in F  \quad \forall 1 \le r \le r, (r, f) \not\in G \label{eqn:rooms1:feature} \\ 
Z_{u,v,r} \ge 0                                    & \quad \forall v \in V  \quad \forall 1 \le f \le f_{\max}
\label{eqn:rooms1:nonnegative}
\end{align}

This results, however, in relaxations too large to be handled by solvers currently available. 
An alternative approach uses fewer additional variables: 
\begin{align}
R_{v,r} = \begin{cases}
\; 1 & \text{ if vertex $v$ is assigned (room) $r$} \\
\; 0 & \text{ otherwise. }
\end{cases}
\end{align}  

The additional constraints follow. In timetabling language: each event is in exactly one room (\ref{eqn:rooms2:iff1}, \ref{eqn:rooms2:iff2}),
events in the same timeslot do not share rooms (\ref{eqn:rooms2:exclusive}),
and $R_{v,r} = 0$ if the event-room combination does not match the event's room-feature or capacity requirements (\ref{eqn:rooms2:capacity},\ref{eqn:rooms2:feature}). 
\begin{align}       
\sum_{1 \le r \le m} R_{v,r} = t   & \quad \forall v \in V \label{eqn:rooms2:iff1} \\
R_{v,r} + R_{v,r'} \le t           & \quad \forall v \in V  \quad \forall 1 \le r, r' \le m, r \neq r' \label{eqn:rooms2:iff2} \\
R_{u,r} + R_{v,r} + Y_{u,v} \le 2t  & \quad \forall u,v \in V, u \neq v  \quad \forall 1 \le r \le m \ \label{eqn:rooms2:exclusive} \\
R_{v,r} = 0                        & \quad \forall v \in V \quad \forall 1 \le r \le m, p_v \ge r_r \label{eqn:rooms2:capacity} \\
R_{v,r} = 0                        & \quad \forall v \in V \quad \forall 1 \le f \le f_{\max}, (v, f) \in F  \quad \forall 1 \le r \le r, (r, f) \not\in G \label{eqn:rooms2:feature} \\
R_{v,r} \ge 0                      & \quad \forall v \in V \quad \forall 1 \le r \le m
\end{align}

Notice that the encoding makes it possible to formulate room stability constraints and penalties,
as it is invariant to ``timeslot permutations''. For example, the hard constraint reads
$R_{v,r} + R_{v',r'} \le t$ for all suitable $v \neq v'$ and all $r \neq r'$.

\section{Algorithms}
\label{sec:algo}

While our key contribution are the actual relaxations, we showcase how these can be used in state-of-the-art algorithms. 
Such algorithmic applications of the relaxations underlie not only our computational results, but also our analytical results in Section \ref{sec:analysis}.

\subsection{Solving the Relaxations}

First, let us consider a first-order method based on the alternating-direction method of multipliers (ADMM) on an augmented Lagrangian, following the extensive literature   
\cite{MR2507127,MR2197554,MR2600237,MR2266705,WenGoldfarbYin2010,GoldfarbMa2010,Yang2015}. 
In order to distinguish between equality constraints reflecting the structure of the conflict graph ($A_1$) and the remainder of the equality constraints ($A_2$), let us consider the primal-dual pair:
\begin{align}
z_p =  \min_{X \in \S^n}  \scal{C}{X} \s.t.\; \OP_{A_1}(X) = b_1 \mbox{ and } \OP_{A_2}(X) = b_2 \mbox{ and } \OP_B(X) \ge d \mbox{ and } X \succeq 0 \notag \\  
z_d = \max_{y_1 \in \R^m, y_2 \in \R^p, v \in \R^q, S \in \S^n} b_1^T y_1 + b_2^T y_2 + d^T v \notag \\
\s.t. \; \OP_{A_1}^*(y_1) + \OP_{A_2}^*(y_2) + \OP_B^*(v) + S = C \mbox{ and } S \succeq 0 \mbox{ and } v \ge 0. \label{eqn:reprinted} 
\end{align}
with the linear operator $\OP_A(X)$ mapping matrix $X$ and matrix $A$ to vector as in the definition of SDPs in Section \ref{sec:sdp-def}.
The augmented Lagrangian of the dual (\ref{eqn:reprinted}) is then:
\begin{align}
L_{\mu}(X, y_1, y_2, v, S) = & - b_1^T y_1 - b_2^T y_2 - d^T v \\
& + \scal{X}{\OP_{A_1}^*(y_1) + \OP_{A_2}^*(y_2) + \OP_B^*(v) + S - C} \notag \\
& + \frac{1}{2\mu} ||\OP_{A_1}^*(y_1) + \OP_{A_2}^*(y_2) + \OP_B^*(v) + S - C||^2_F. \notag
\end{align} 
In an alternating direction method of multipliers, one minimises the augmented Lagrangian in $v, S$, and $(y_1, y_2)$, in turns, as suggested in Algorithm Schema~\ref{algo:bpm}.

\begin{algorithm}[t!]
\caption{{\tt AugmentedLagrangianMethod}($A_1, A_2, B, C, b_1, b_2, d$) } 
\label{algo:bpm}
\begin{algorithmic}[1]
  \STATE \textbf{Input:} Instance I = ($A_1, A_2, B, C, b_1, b_2, d$) of SDP, precision $\epsilon$(\ref{eqn:reprinted})  
  \STATE \textbf{Output:} Primal solution $Y$, computed up to $\epsilon$-precision \rule{0pt}{2.5ex}
  \vspace{1ex}
\STATE Set iteration counter $k = 0$ 
\STATE Initialise $X^{k} \succeq 0$ with a heuristically obtained colouring
\STATE Compute matching values of dual variables $y_1^{k}, y_2^{k}, v^{k} \ge 0$, and $S^{k} \succeq 0$ 
\WHILE{ the precision is insufficient  } \label{line:loop}
\STATE \begin{tabular}{l} Increase iteration counter $k$ \label{line:inc} \end{tabular}
\STATE \begin{tabular}{ll}
Update $v^{k+1}$ & = $\argmin_{v \in \R^q, v \ge 0} L_{\mu}(X^{k},  y_1^{k+1}, y_2^{k+1}, v, S^{k})$
\end{tabular}                 
       \label{line:QP}
\STATE \begin{tabular}{ll} 
Update $S^{k+1}$ & = $\argmin_{S \succeq 0} L_{\mu}(X^{k}, y_1^{k+1}, y_2^{k+1}, v^{k+1}, S)$ 
\end{tabular} \label{line:eig}                                  
\STATE \label{line:closed-form}
\begin{tabular}{ll}
$(y_1^{k+1}, y_2^{k+1})$ & = $\argmin_{y_1 \in \R^m, y_2 \in \R^m} L_{\mu}(X^{k}, y_1, y_2^{k}, v^{k}, S^{k})$
\end{tabular} 
\STATE \begin{tabular}{l} Choose any step-length $\mu \ge 0$ \end{tabular}  
\STATE \begin{tabular}{ll} Update $X^{k+1}$ & = $X^{k} + \frac{A_1^T(y_1^{k+1}) + A_2^T(y_2^{k+1}) + B^T(v^{k+1}) + S^{k+1} - C)}{\mu}$ \end{tabular} 
\ENDWHILE
\STATE Return $X$
\end{algorithmic}
\end{algorithm} 

In general, Algorithm Schema~\ref{algo:bpm} reduces the minimisation of one moderately 
complicated convex optimisation problem to solving three simpler convex optimisation sub-problems.
In Line \ref{line:QP}, one can solve the linear system given by first-order Karush–Kuhn–Tucker optimality conditions of 
\begin{align}					
\argmin_{v \in \R^q, v \ge 0} \left( \left( {B\left( {X^k  + \frac{1}{{\mu}}\left( {A_1^T (y_1^{k + 1} ) + A_2^T (y_2^{k + 1} ) + S^k - C} \right)} \right) - d} \right)^T v + \frac{1}{{2\mu}}v^T (BB^T ) v \right).
\end{align}
In Line \ref{line:eig},
it is important to realise that 
\begin{align}
\argmin_{S \in \S^n, S \succeq 0 } \norm{S - \left(C - A_1^T(y_1^{k+1}) - A_2^T(y_2^{k+1}) - B^T(v^{k+1}) - \mu X^k \right)}_F^2
\end{align}
can be solved by spectral decomposition of the term subtracted from $S$ \cite{MR1247916}.
Finally, in Line \ref{line:closed-form}, one can initialise the computation with:
\begin{align}
y_1^{k+1} = & -(A_1 A_1^T)^{-1} (\mu(A_1(X^k) - b_1) + A_1(A_2^T(y_2^k) + B^T(v^k) + S^k - C))\\
y_2^{k+1} = & -(A_2 A_2^T)^{-1} (\mu(A_2(X^k) - b_2) + A_2(A_1^T(y_1^{k+1}) + B^T(v^k) + S^k -C)).
\end{align}
We refer to \cite{Yang2015} for a some excellent suggestions as to the implementation of the linear solver and spectral decomposition, as well as convergence properties
of such as method.

\begin{algorithm}[t!]
\caption{{\tt Rounding}($X$) based on Karger, Motwani, and Sudan}
\label{algo:rounding}
\begin{algorithmic}[1]
  \STATE \textbf{Input:} Matrix variable $X$ of the solution to the SDP (\ref{OUR-SDP}) of dimensions $n \times n$,
                         bound $m$, number $a_{\max}$ of randomisations to test,
                         plus the input to Simple Timetabling, if required 
  \STATE \textbf{Output:} Partition $P$ of the set $V = {1, 2, \ldots, n}$ \rule{0pt}{2.5ex}
  \vspace{1ex}
\STATE Compute vector $v, X = v^T v$ using Cholesky decomposition
\FOR{ Each attempted randomisation $a = 1, \ldots, a_{\max}$}
\STATE Initialise $P_a = \emptyset, i=1, X = V$  
\WHILE{ There are uncoloured vertices in $X$ }
\STATE Pick a suitable $c = \sqrt{ \frac{2(k - 2)}{k \log_e \Delta} }$ for $\Delta$ being the maximum degree of the vertices in $X$ 
\STATE Generate a random vector $r$ of dimension $|X|$ 
\STATE Pick $R_i \subseteq X$ of at most $m$ elements in the descending order of $v_i r_i$, 
          where (1) positive and (2) independent of previously chosen and, in Simple Timetabling, 
          (3) the respective events fit within the rooms and (4) require only features available 
\STATE Update $P_a = P_a \cup \{ \{ R_i \} \}, X = X \setminus R_i$, $i = i + 1$
\ENDWHILE
\ENDFOR
\STATE Return $P_a$ of minimum cardinality 
\end{algorithmic}
\end{algorithm}

\subsection{Recovering an Assignment}
\label{sec:recoveryAlg}

Let us comment on the recovery of an upper bound from the lower bound provided by SDP.
Since the seminal paper of Karger, Motwani, and Sudan \cite{Karger1998}, 
there has been a continuing interest in algorithms recovering a colouring from
semidefinite relaxations. 
Typically, such algorithms are based on simple randomised 
iterative rounding of the semidefinite programming relaxation. One such algorithm, specialised to simple timetabling is displayed in Algorithm Schema~\ref{algo:rounding}. 

Alternatively, one can consider methods solving a sequence of smaller semidefinite programming relaxations, inspired by the so-called iterated rounding in linear programming \cite{lau2011iterative}. 
When applied to linear programming, the method fixes variables, whose values in the relaxation are close to $0$ or $1$, to 0 or 1, respectively, and resolves the smaller residual linear program.
When applied to semidefinite programming, the method 
fixes eigenvectors whose corresponding eigenvalues are close to zero or one. 
Let us consider the example of \cite{morgenstern2019fair} starting from:
\begin{align}
  & \min \langle C, X\rangle  & \label{eq:suitable1} \\
 \st \quad &  \langle A_i , X\rangle  \geq   b_i & \;\; \forall \; 1\leq i\leq m \notag \\
   & \trace(X) \leq d & \notag \\
   & 0\preceq X \preceq  I_n, & \notag 
\end{align}
which can accommodate many of the SDP relaxations we have seen so far.
There, \cite{morgenstern2019fair} initialise  $F_0=F_1=\emptyset$ and $F=I_n$. 
In each iteration, subspaces spanned by eigenvectors corresponding to eigenvalues $0$ or $1$ are fixed and the corresponding standard basis vectors are moved from $F$ to $F_0$ and $F_1$, respectively. 
Thus, one increases the subspaces spanned by columns of $F_0$ and $F_1$, while maintaining pairwise orthogonality. 
To obtain new $F$, one solves a smaller semidefinite program in $r\times r$ symmetric matrix $X(r)$:
\begin{align}
\label{eq:smallersdp}
\max \, & \, \langle F^TCF, X(r)\rangle \\
\langle F^T A_i F, X(r)\rangle & \ge b_i- F_1^T A_i F_1 \quad i \in S \notag \\
\trace(X(r)) & \le d-\rank(F_1) \notag \\
0 \preceq \, X(r)  &\preceq I_r, \notag
\end{align}
which assures that, eventually, we can recover $X$ that is orthogonal to all vectors in subspace spanned by vectors in $F_0$,
and whose eigenvectors corresponding to eigenvalue $1$
will be the columns of $F_1$.
This is summarised in Algorithm Schema~\ref{algo:rounding2}. 
As we will see in Section~\ref{sec:recovery}, this allows for non-trivial performance guarantees.

\begin{algorithm}[t!]
\caption{{\tt IterativeRounding}($X$) based on Morgenstern et al.}
\label{algo:rounding2}
\begin{algorithmic}[1]
  \STATE \textbf{Input:} An $n \times n$ matrix  $X$ of the solution to the SDP (\ref{eq:suitable1}),
 which has $m$ inequalities, alongside with the corresponding matrices $A_i$ for $i = 1, \ldots, m$ 
  \STATE \textbf{Output:} Partition $P$ of the set $V = {1, 2, \ldots, n}$ \rule{0pt}{2.5ex}
  \vspace{1ex}
\STATE Initialize $F_0, F_1$ to be empty matrices and $F=I_n$, $S\gets \{1,\ldots, m\}$.
\STATE Initialise $\delta > 0$ to be a threshold for rounding
\WHILE{ $F$ is non-empty }
\STATE Solve \eqref{eq:smallersdp} to obtain extreme point $X^*(r)=\sum_{j=1}^r \lambda_j v_j v_j^T$ where $\lambda_j$ are the eigenvalues and $v_j\in \R^r$ are the corresponding eigenvectors.
\STATE For any eigenvector $v$ of $X^*(r)$ with eigenvalue less than $\delta$, let $F_0\gets F_0\cup \{Fv\}.$
\STATE For any eigenvector $v$ of $X^*(r)$ with eigenvalue of more than $1 - \delta$, let $F_1\gets F_1\cup \{Fv\}.$
\STATE Let $X_f=\sum_{j: 0<\lambda_j<1} \lambda_j v_j v_j^T$. If there exists a constraint $i\in S$ such that $\langle F^T A_i F, X_f\rangle < \delta$, then
$S\gets S\setminus \{i\}.$
\STATE Update $F$ by taking every eigenvector $v$ of $X^*(r)$ with eigenvalue within $[\delta, 1-\delta]$, and taking $Fv$ to be the columns of $F$.
\ENDWHILE
\STATE From rank-$t$ matrix $F_1F_1^T$ reconstruct partition $P$ by Cholesky decomposition
\end{algorithmic}
\end{algorithm}

As a remark, we note that there are many other alternative rounding approaches within the Theoretical Computer Science literature. We refer to \cite{6108208,Raghavendra2012,Bansal2016,abbasi2018sticky} for notable examples. While they may not be directly applicable, they are based on important insights that would be applicable.

\section{An Analysis}
\label{sec:analysis}

Next, let us analyse the strength of the bound and the complexity of computing it, both of which affect its practicality.

\subsection{The Strength of the Bound}

In terms of strength of the bound, one can extend a number of properties of relaxations of graph colouring to bounded colouring. For the sake of completeness, we reiterate some of them. 
For instance, one can show the sandwich-like:
\begin{prop}
\label{prop:sandwich}
For every graph $G$, there is an $m \ge 0$, such that 
\begin{align}
\omega(G)   \le \LB'(\overline G)   
            \le \chi(G)
            \le \mathcal{C}(G, m)
            \le \OUR'(G, m)
            \le \OUR(G, m) \\
			\OUR(G, m)
            \le \OUR^+(G, m)
            \le \OUR^{+\bigtriangleup}(G, m)
            \le \chi(G, m),
\end{align}
\end{prop}
where $\omega$ is the size of the largest clique,
$\chi$ is the chromatic number,
$\mathcal C$ is a bound obtained by counting,
$\chi^m$ is the $m$-bounded chromatic number,
the values of SDP relaxations follow the notation of Figure~\ref{over1} and
$\OUR^{+\bigtriangleup}$ is the strengthening of $\OUR^+(G,m)$ with triangle inequalities.
\begin{sketch}
The relationship between the values of the successive relaxations of bounded colouring is clear.
To show there is $m$, such that $\chi(G) \le \mathcal C(G, m)$, let us study two cases:
If there is $r \ge 1$ such that $r$-bounded
graph colouring of $G$ requires a strictly larger number of colour classes than the chromatic number, take $m = r$.
Otherwise, the graph cannot have independent sets larger than one, hence is a clique, and $\omega(G) = \chi(G) = \chi(G, m)$ for any $m$.
\end{sketch}

To see that (non-bounded) graph colouring relaxations   
($\LB(\overline G), \LB'(\overline G), \LB^+(\overline G)$)
provide only very weak bounded graph colouring relaxations,
consider empty graphs on $n$ vertices and the constant function $f(n) = 1$:  
\begin{prop}
There is an infinite family of graphs and $f(n)$, where the chromatic number is $O(1)$,
the $f(n)$-bounded chromatic number is $O(n)$.
\end{prop}
In contrast, the value of the semidefinite programming relaxation of bounded colouring 
may match the bounded chromatic number on such graphs.

On random graphs where an edge between each pair of distinct vertices appears with probability $p$, independent of any other edge, which are known as Erd\H{o}s-R{\' e}nyi $G(n,p)$:

\begin{prop}
\label{prop:Juhasz}
With probability $1 - o(n)$, graph $G$ drawn randomly from $G_{n,p}$ has
\begin{align}
\OUR(G) \ge \frac{\sqrt{n} }{2} \sqrt{\frac{1 - p}{p}} + O(n^{\frac{1}{3}} \log n),
\end{align} 
where the big-$O$ notation hides lower-order terms. 
\end{prop}

\begin{sketch}
The proof combines the sandwich-like Property~\ref{prop:sandwich} and the impressive result of \citeA{MR685042}.
\end{sketch}

Computationally, this bound seems to be rather tight, as we show in Section \ref{sec:computational}.

\subsection{The Structure of the Relaxations}
\label{sec:structure}

For example, let us consider the formulation of bounded graph colouring (\ref{eqn:Rewritten}) for a graph on $n$ vertices 
and $m$ edges. There equality constraints reflecting the structure of the conflict graph ($A_1$) have the cardinality of their support (number of non-zero elements) equal to the number of edges in the conflict graph and the remainder of the equality constraints ($A_2$) also have a very simple structure:

\begin{prop}
\label{prop:a1a1t}
First $m$ equalities (\ref{eqn:A1}) correspond to $m \times n^2$ matrix $A_1$.
$A_1 A_1^T = I_m$, where $I_m$ is the $m \times m$ identity matrix.
\end{prop}

\begin{prop}
\label{prop:a2a2t}
Further $n-1$ equalities (\ref{eqn:A2}) correspond to $n-1 \times n^2$ matrix $A_2$.  
$A_2 A_2^T = J_{n-1} + I_{n-1}$, where $I_{n-1}$ and $J_{n-1}$ are $(n-1) \times (n-1)$ identity and all-ones matrices, respectively.
$(A_2 A_2^T)^{-1}$ is $- \frac{1}{n} J_{n-1} + I_{n-1}$.
For $(n-1)$-vector $y$, $A_2^T y$ is 
an $n \times n$ matrix, with $\left[\left(\sum_i y_i\right) (-y_1) (-y_2) \cdots (-y_{n-1})\right]$ on the diagonal and zeros elsewhere.  
For positive $X$, $\OP_{A_2}(X) = \left[ 2, 4, \cdots, 2(n-1) \right]$ of dimension $(n - 1)$. 
\end{prop}

\begin{prop}
\label{prop:btv}
Inequalities (\ref{eqn:B}) correspond to $n \times n^2$ matrix $B = I \kron j$, 
where $j$ is the row-vector of $n$ ones. Hence, $BB^T = n I_n$,
where $I_n$ is the $n \times n$ identity matrix.
For an $n$-element column-vector $v$, $B^T v = (v \kron j)^T = \left[v_1 j \; v_2 j \; \cdots \; v_n j \right]^T$,
where $j$ is the row-vector of $n$ ones. 
\end{prop}

\begin{prop}
The elements of the objective matrix $C$ are zeros except for $C_{1,1} = 1$. 
Hence $\OP_{A_1}(C) = 0$, where $0$ is the $m$-vector of zeros.
$\OP_{A_2}(C) = j$, where $j$ is the $(n - 1)$-vector of ones.
\end{prop}

Across both:
\begin{itemize} 
\item custom solvers, such as Algorithm Schema \ref{algo:bpm} and the three sub-problems in Lines \ref{line:QP}--\ref{line:closed-form}, in particular, and 
\item general-purpose solvers allowing for the input of block-structured matrices with sparse and identity blocks, such as \cite{Fujisawa2000,Gondzio2009}, 
\end{itemize}
it is possible to exploit Properties~\ref{prop:a1a1t}--\ref{prop:btv} so as to:  
\begin{itemize}
\item not compute $(A_1 A_1^T)^{-1}$
\item compute $A_1^T y_1$ in time $O(m)$
\item compute $(A_2 A_2^T)^{-1}$ in time $O(n^2)$
\item compute $A_2^T y_2$ in time $O(n)$
\item compute $(A B^T)^{-1}$ in time $O(n)$  
\item compute $B^T v$ in time $O(n)$ 
\item evaluate the augmented Lagrangian and its gradient at a given $v$ in time $n^2$    
\end{itemize}
in relaxations of bounded graph colouring of a graph on $n$ vertices,
compared to $O(n^6)$ run-time of methods not exploiting the structure.

\subsection{The Recovery}
\label{sec:recovery}

Due to the hardness of approximation of colouring in a graph with large enough a chromatic 
number within the factor of $n^\epsilon$ for some fixed $\epsilon$ \cite{MR2403018}, 
one cannot hope to guarantee reconstruction of a solution close to optimality in the worst case.
Having said that, as we will illustrate in the next section, however, Algorithm Schema~\ref{algo:rounding} performs rather well in practice. 

One can also provide weaker guarantees.
In particular, one could consider the so-called \(\epsilon\)-solution, 
 which satisfies linear constraints within an additive error of $\epsilon$, while being at most $\epsilon$ from the optimal objective.
Notice that the fact that an \(\epsilon\)-solution is obtainable in 
time polynomial in $n$ and $\log \frac{1}{\epsilon}$
does not contradict the hardness of approximation results \cite{MR2403018}, which consider the objective of solutions satisfying the constraints exactly. 

\begin{prop} 
There exists an $\epsilon > 0$  and an algorithm implementing Algorithm Schema \ref{algo:rounding2} 
that, given any feasible solution 
to the SDP relaxation \eqref{OUR-SDP-WithFeatures} of {\sc Simple Laminar Timetabling},
runs in time 
polynomial in $n$ and $\log \frac{1}{\epsilon}$
and returns an $\epsilon$-feasible and \(\epsilon\)-optimal solution to the SDP  relaxation \eqref{OUR-SDP-WithFeatures} of {\sc Simple Laminar Timetabling}.
\end{prop}

\begin{sketch}
The proof extends the work of  \cite{morgenstern2019fair} on the number of fractional eigenvalues in any extreme point  $X$ of a suitable form of a semidefinite program with $m$ linear  inequalities and trace bounded by $t$, which is
\begin{align}
    t + \floor*{ \sqrt{2m+\frac{9}{4}}-\frac{3}{2} 
    }.
\end{align}
Based on this bound, one can formulate a generic result on iterative rounding of SDPs, which we present in Proposition \ref{thm:SDProunding} below.
The result applies to the 
SDP relaxation \eqref{OUR-SDP-WithFeatures} of {\sc Simple Laminar Timetabling}, because it can be cast into the suitable form \eqref{eq:suitable1}. 
The bound on the run-time follows from the fact we solve at most $n$ semidefinite programs in matrices at most $n \times n$ and standard results on interior-point methods \cite{MR1315703}.
\end{sketch}

\begin{prop}[Theorem 7 of \cite{morgenstern2019fair}]
\label{thm:SDProunding}
Let $C$ be a $n\times n$ matrix and $\{A_1,\ldots, A_m\}$ be a collection of $n \times n$ real matrices, $d \leq n$, and $b_1,\ldots, b_m\in\R$. Suppose the semi-definite program \eqref{eq:suitable1}
with a trace bounded by $d$ and $m$ other constraints has a nonempty feasible set and let $X^*$ denote an optimal solution. There is an algorithm 
  that given a matrix $X_0$ that is a strictly feasible solution, 
  returns a matrix $\tilde{X}$ such that
 \begin{enumerate}
 \item rank of $\tilde{X}$ is at most $d$,
 \item  $\langle C, \tilde{X}\rangle \leq \langle C, X^* \rangle $, and
  \item for each index $1\leq i\leq m$ of a constraint we have \begin{align}
  \label{eq:violation-bound}
  \langle A_i , \tilde{X}\rangle   \geq   b_i- \max_{S\subseteq [m]} \sum_{i=1}^{\floor*{ \sqrt{2|S|}+1 }} \sigma_i(S),\end{align}
  where $\sigma_i(S)$ is the $i^{th}$ largest singular of the average of matrices $\frac{1}{|S|} \sum_{i\in S} A_i$ for any subset of matrices defining the constraints,  $S\subseteq \{1,\ldots, m\}$.
 \end{enumerate}
\end{prop}

In the violation bound \eqref{eq:violation-bound}, the quantity $\sigma_i(S)$ is  non-trivial to reason about, but it is clear that it is rather modest, because the singular values are at most 1 and the summation goes over at most $\sqrt{2m}+1$ values.

\section{Computational Experience}
\label{sec:computational}

To corroborate our analytical results in Section \ref{sec:analysis}, we have conducted a variety of computational tests. 
Most of these have been driven by YALMIP \cite{YALMIP} 
scripts running within MathWorks Matlab R2017b on a laptop with 
Intel Core Duo i5 at 2.7~GHz with 8~GB of RAM,
which also had 
IBM ILOG CPLEX 12.8 and
SeDuMi 1.3 \cite{SEDUMI} installed. Let us refer to it as a laptop.
When explicitly mentioned, we also present results
obtained on a machine equipped 
with 80 cores of Intel Xeon E7-8850 at 2.00~GHz
and 700~GB of RAM, which had
MathWorks Matlab R2016b, 
IBM ILOG CPLEX 12.6.1,
and SeDuMi 1.3 \cite{SEDUMI}
installed.
Let us refer to it as a large-memory machine.

\subsection{A Motivating Example}

As a first concrete computational example, we consider a small conflict graph from a standard collection of benchmark problems in timetabling. Specifically, we take the instance {\tt sta-f-83} from the Toronto examination timetabling benchmarks \footnote{See \texttt{ftp://ftp.mie.utoronto.ca/pub/carter/testprob/} and \texttt{http://www.cs.nott.ac.uk/$\sim$rxq/data.htm}}. There are 139 events, but the conflict graph has three connected components of 30, 47 and 62 vertices. Here, we use the 47-vertex component. 
The results are given in Table~\ref{tab:bounded},
with bounded chromatic numbers obtained using the most straightforward integer linear programming formulation 
solved using the default settings of IBM ILOG CPLEX on a laptop.

\begin{table}[t!]
%\begin{center}
\caption{
An illustration of the effects of bounding the $m$-bounded chromatic number of the instance {\tt sta-f-83}: 
Column $\chi^{m}$ lists the $m$-bounded chromatic number obtained using integer linear programming, within time listed under ``$\chi^{m}$ Runtime'' in seconds.
Column $\OUR^{m}$ lists the bounds obtained using semidefinite programming and rounding up, within time listed under ``$\OUR^{m}$ Runtime'' in seconds.
% *** TODO *** DO NOT ROUND UP
Column $|V|/m$ lists the lower bound on the colours obtained by simple counting arguments and rounding up.  
Dash denotes the the omission of the $m$-bounding constraint, giving $\LB$ instead of $\OUR^{m}$.
}
\label{tab:bounded}
  \begin{tabular}{l|ll|ll|l}
$m$ & 
$\chi^{m}$   &  $\chi^{m}$ Runtime &
$\OUR^{m}$   &  $\OUR^{m}$ Runtime &
$|V|/m$ \\
\hline
1 & 47 & 0.09 & 47 & 3.46 &47 \\ 
2 & 26 & 2.88 & 26 & 2.92 & 24\\ 
3 & 20 & 2.67 & 20 & 3.34 & 16 \\ 
4 & 16 & 7.22 & 16 & 3.70 & 12 \\ 
5 & 14 & 11.10 & 14 & 3.24 & 10 \\ 
6 & 13 & 2.67 & 13 & 3.12 & 8 \\ 
7 & 12 & 8.77 & 12 & 3.26 & 7 \\ 
8 & 11 & 2.89 & 11 & 3.40 & 6\\ 
9 & 11 & 3.39 & 11 & 3.14 & 6\\ 
47 & 11 & 0.35 & 11 & 3.92 & 1 \\ \hline 
--- & 11 & 0.34 & 11 & 3.45 & --- \\ \hline 
\end{tabular}
%\end{center}
\end{table}

Firstly, note that $m=1$ gives precisely the number of nodes, as would be expected. 
Secondly, note that $\OUR^m$ is generally much tighter than the lower bound $|V|/m$ obtained by simple counting arguments.
Accidentally, $\OUR^m$ lower bounds actually happen to match the optima in this particular instance. 
For example, at $m=5$, counting cannot rule out a 10-colouring, but the SDP bound shows that at least 14 colours are required. 
As far as we know, SDP relaxations are the only way to get such information in polynomial time, considering that
the 14-colouring together with a certificate of its optimality can be obtained using CPLEX, but not in polynomial time.

\subsection{Random Graphs}

Next, we show that the same behaviour can be observed on a large sample of random graphs. 

\begin{figure}[t!]
\caption{The effects of tightening the bound on the number of uses of a colour on the strength of the lower bound:
For a random graph $G(n, 0.5)$, where the size of the largest colour class in an optimal colouring is $C$, 
the mean lower bounds, upper bounds, and optima for 
unbounded colouring, 
$(C-1)$-bounded colouring, $(C-2)$-bounded colouring, etc.,
are computed from a sample of $N = 100$ for each number of vertices $n$ and restriction on the size of the colour class.}
\label{fig:random1}
\centering
\includegraphics[width=0.9\textwidth]{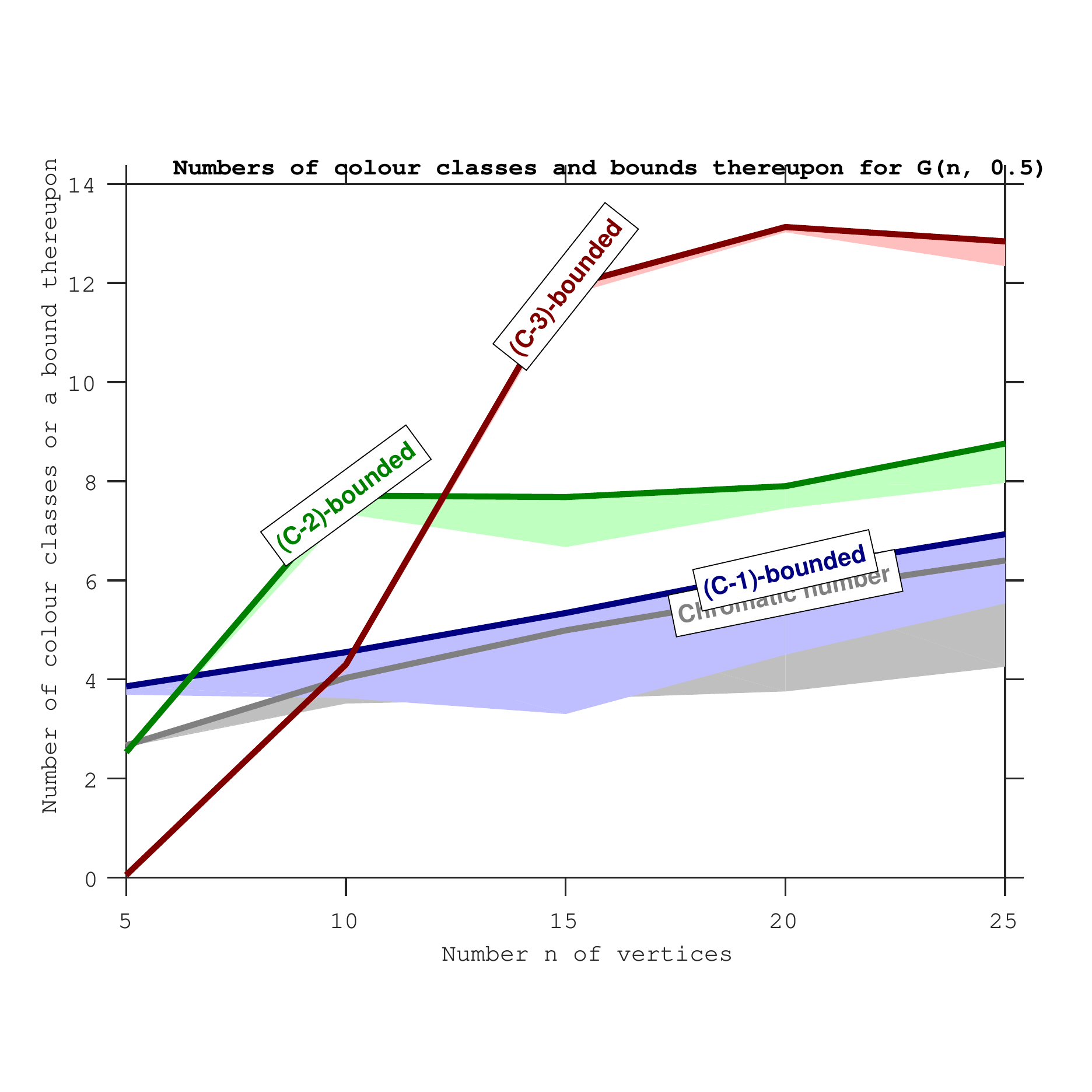}
% random-plots-n4-pretty-solid
\end{figure}

First, we demonstrate the improved strength of the lower bound obtained from semidefinite programming
 as the restriction on the number of uses of a colour is tightened (i.e., cardinality of a colour class is bounded from above by progressively smaller numbers). 
In general, we compute the best possible vertex colouring, without any bound on the number of uses of a colour,
and take the size of the largest colour class to be $C$. 
Subsequently, we obtain lower bounds, upper bounds, and optima for $(C-1)$-bounded colouring, $(C-2)$-bounded colouring, etc., of the same graph.
In particular, we use random graphs with constant probability 0.5 of an edge appearing between a pair of distinct vertices and varying numbers $n$ of vertices, which are known as $G(n, \frac{1}{2})$. 
For each number $n$ of vertices, we have generated 100 random graphs, computed the true chromatic numbers and 
the size of the largest colour class 
$C$ using CPLEX,  
lower bounds on the bounded colouring using SeDuMi, 
and upper bounds by rounding the semidefinite programming relaxation, all running on a laptop.
%(The data are available\footnote{
%Please see \url{http://github.com/***/} (Retrieved January 5th, 2018).
%} on-line.)
%Notice the fit with the prediction of Juhasz (\ref{prop:Juhasz}) in the
%the means over the 100 sample paths, which are presented in 
In Figure~\ref{fig:random1},
the true value is plotted in a solid line, while a semi-transparent region spans the lower and upper bounds.
Notice that:
In Figure~\ref{fig:random1},
the true value is plotted in a solid line, while a semi-transparent region spans the lower and upper bounds.
Notice that:
\begin{itemize}
\item the upper bounds obtained by rounding the semidefinite programming relaxation coincides with the true value obtained by ILOG 
\item for unbounded and $(C-1)$-bounded colouring, there is a considerable gap between the SDP-based lower bound and the true value
\item for $(C-3)$-bounded colouring, the SDP-based lower bound and the true value coincide in the majority of cases. (That is: One can round the upper bound up, as it has to be integral. The average over 100 samples need not be integral, though.)
\item for $(C-3)$-bounded colouring, the SDP-based lower bound is essentially tight. 
\end{itemize}

\begin{figure}[t!]
\caption{The run-time of the presented methods on a laptop as a function of the number of vertices:
Sample mean run-times of 
an interior point method (IPM) and an 
augmented Lagrangian (AugLag) method on relaxations for $G(n, 0.5)$ on a laptop
for $N = 100$ samples per each number $n$ of vertices.}
\label{fig:random3}
\centering
\includegraphics[width=0.9\textwidth]{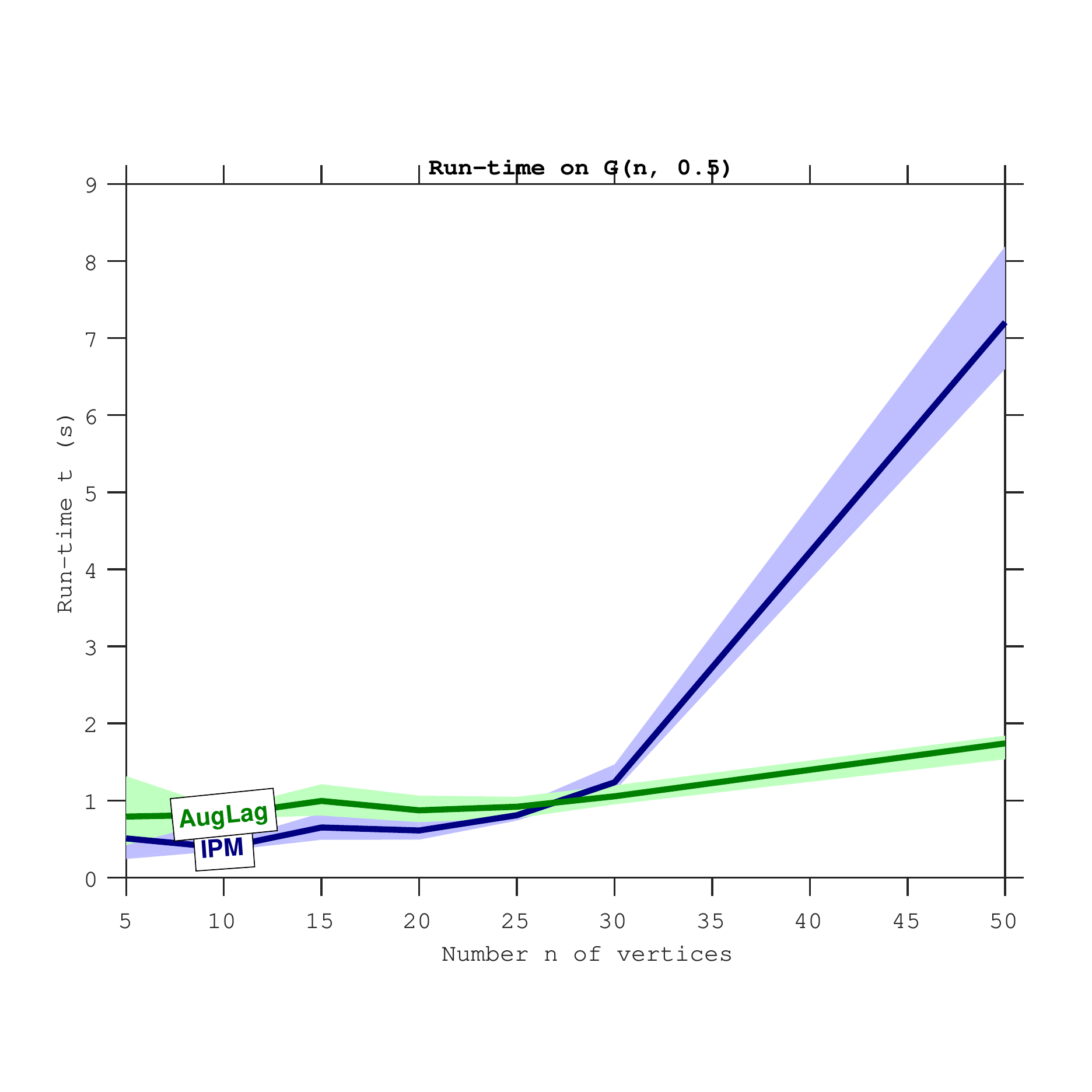}
\end{figure}

Second, we illustrate the practicality of the approach by illustrating the dependence of the run-time on dimensions of the graph.  
Figure~\ref{fig:random3} presents the results on instances, which fit within the memory of a laptop.
It suggests that run-time 
of commonly used
first-order methods for solving semidefinite-programming relaxations increases linearly with the number of vertices of the graph, while the run-time of commonly used 
second-order methods for solving semidefinite-programming relaxations increases quadratically  with the number of vertices of the graph.
This is surprising. Consider the fact that the dimension of the matrix variable increases quadratically with the number of vertices and the number of elements in the Hessian matrix considered in second-order methods increases quadratically in the dimension.
In both cases, the observed run-time is due to the ability of the respective methods to exploit the structure of Section~\ref{sec:structure}.

\begin{figure}[t!]
\caption{The run-time of the augmented Lagrangian (AugLag) method on a large-memory machine, as a function of the number of vertices $n$ in $G(n, 0.5)$. We restrict ourselves to $N = 1$ sample per each number $n$ of vertices, due to the run-time of YALMIP constructing the SDP instances.}
\label{fig:random4}
\centering
\includegraphics[width=0.9\textwidth]{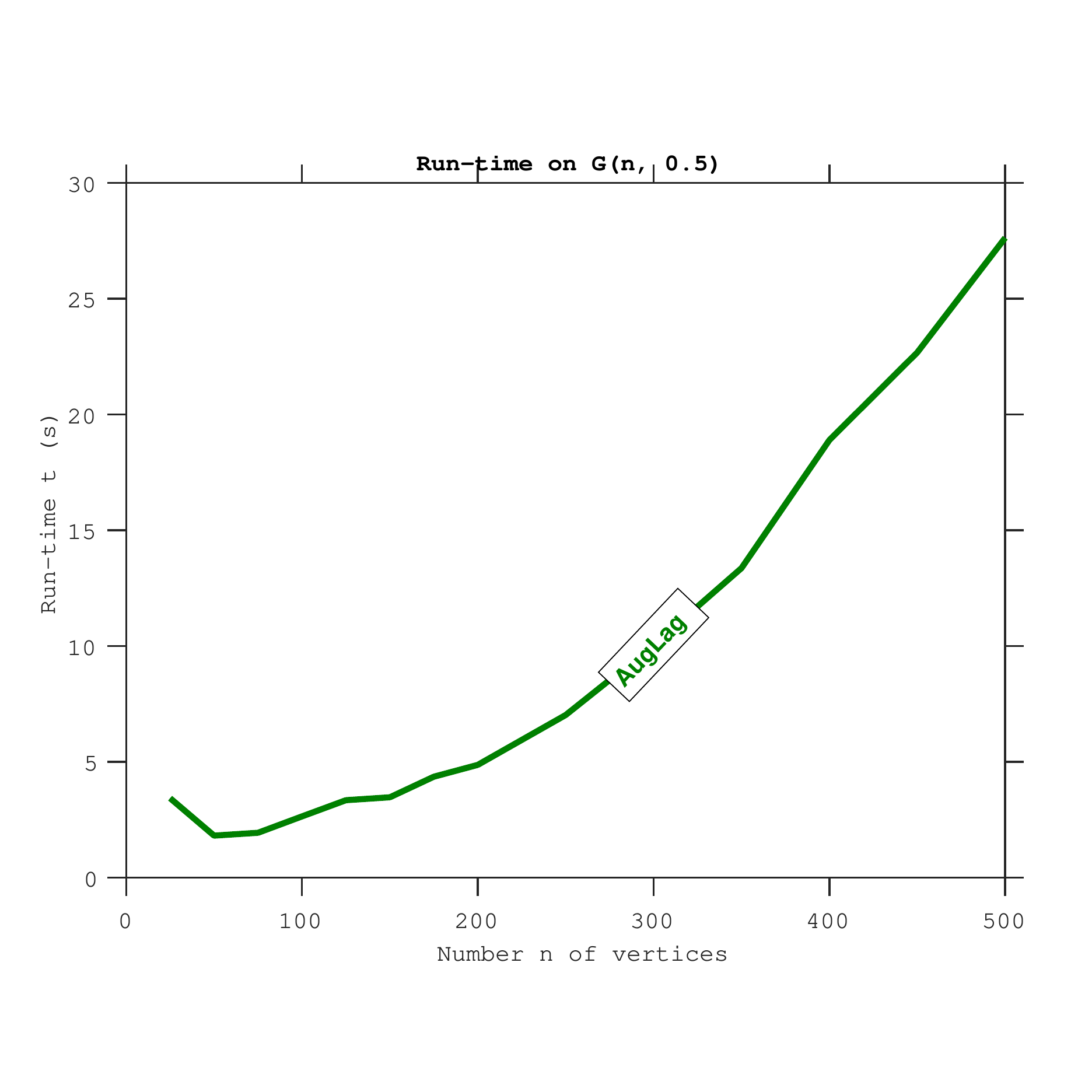}
\end{figure}

Figure~\ref{fig:random4} presents the corresponding results on instances, which no longer fit within the memory of a laptop,
as run on a large-memory machine.
We note that already on a random graph on 200 vertices, $G(200, 0.5)$, SeDuMi 
regularly consumes over 24 GB physical memory, with further 
13 GB in swap,
in solving the SDP relaxation of bounded graph colouring.
Although the run-times are longer, considering the sheer amounts of data processed, the evolution of run-time as a function of the number of vertices seems similar to Figure~\ref{fig:random3}.

\begin{figure}[t!]
\caption{The run-time of the presented methods on a laptop, as a function of graph's density:
For random graphs $G(25, p)$, sample mean run-times of 
an interior point (IPM), possibly with a
with dimension reduction (DimRed) and exploitation of sparsity (SparseCoLo), compared against the run-times of an
augmented Lagrangian (AugLag) method,
for $N = 100$ samples per each density $p = 0.1, 0.2, \ldots, 0.9$.}
\label{fig:random2}
\centering
\includegraphics[width=0.9\textwidth]{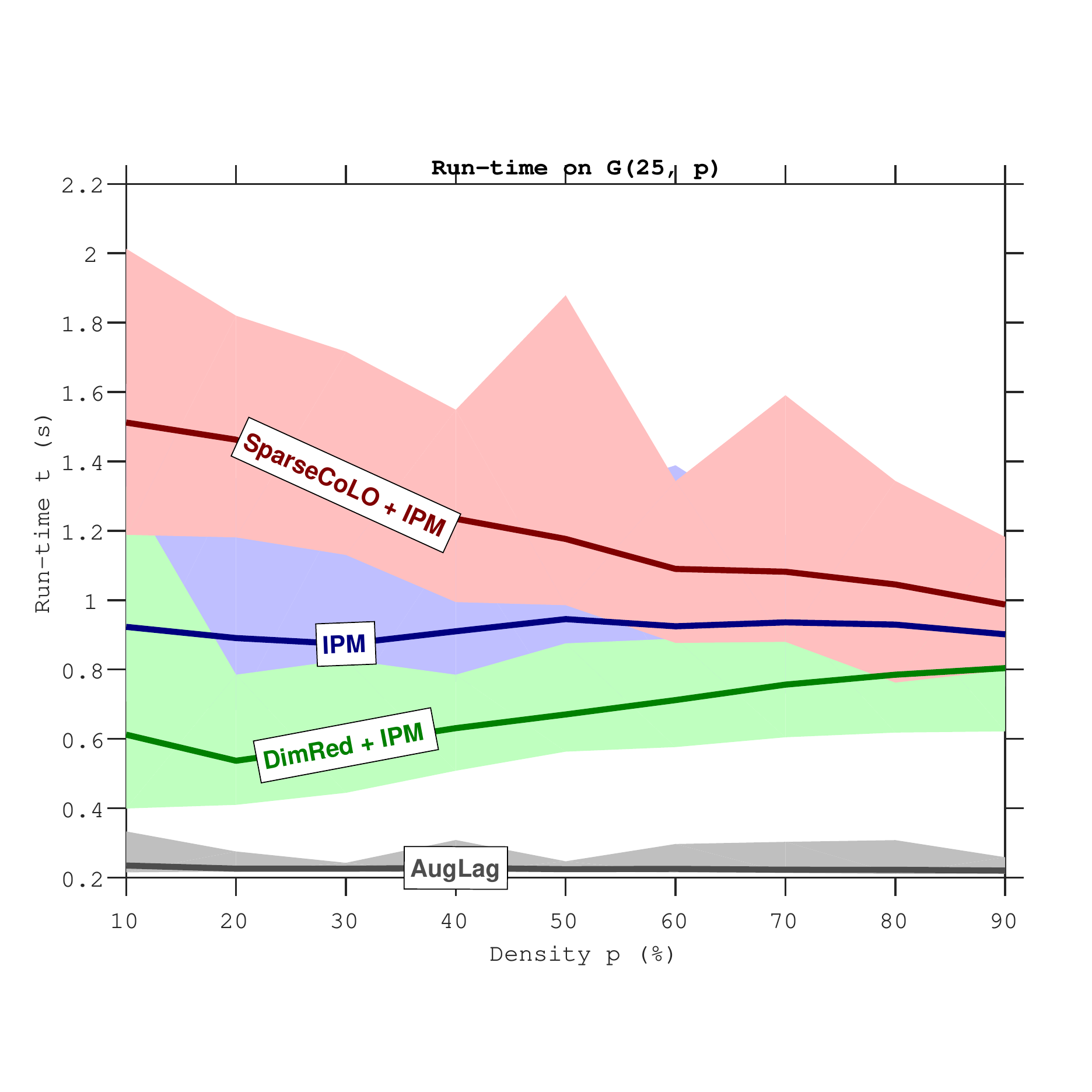} 
% random-runtimes-n40-pretty-solid.pdf
\end{figure}

Figure~\ref{fig:random2} illustrates that commonly used
methods do not exhibit a major increase in run-time as the density of the graph increases, due to their ability to exploit the structure of Section~\ref{sec:structure}. Again, this is surprising. Consider that the number of edges in the conflict graph asymptotically approaches the square of the number of vertices in a dense graph. If there were no structure, the cubic increase of run-time with each of the quadratic number of constraints may render the approach impractical.
In particular, we use random graphs
of varying densities, all on 100 vertices. The data are again available on-line.
For each of the densities $p = 0.1, 0.2, $\ldots$, 0.9 \%$, we have generated 100
 graphs $G(40, p)$.   
Subsequently, we have obtained SDP-based lower bounds using 4 different methods:
\begin{itemize}
\item[IPM]
 a standard implementation of a primal-dual interior point method (IPM) by \cite{SEDUMI}
\item[DimRed] a dimension-reduction procedure of \cite{YALMIP} followed by the IPM of \cite{SEDUMI}
\item[SparseCoLo] a sparsity-exploiting procedure of \cite{kim2011} followed by the IPM of \cite{SEDUMI}
\item[AugLag]
 an implementation of a first-order method considering the augmented Lagrangian, based on \cite{MR2600237,Yang2015}
\end{itemize}
For each of the methods, and each of the densities $p$, we report the average run-time over 100 graphs.

%\begin{table}[t]
%\caption{For event-based conflict graphs from 
%Track 3 of International Timetabling Competition 2007 and Toronto Examination Timetabling Benchmark, where the size of the largest colour class in an optimal colouring is $C$, lower bounds $\OUR^{m}$ and optima $\chi^{m}$ for $(C - m)$-bounded colouring are shown. For $m = 0$, no bounds were applied.}
%\label{tab:Udine}
%\begin{center}
%\include{gfx/Udine-overview}
%\end{table}

\begin{table}[t!]
\caption{Results for instances from %Track 1 (exam), Track 2 (post), and 
Track 3 (comp) of the International Timetabling Competition 2007.
%Notice that unlike in Table~\ref{tab:Udine}, conflict graphs generated from Track 3 are course-based.
}
\label{tab:ITC}
%\begin{center}
\begin{tabular}{l|rr|rr|r}
Graph & $\LB$ & Runtime & $\LB^m$ & Runtime & Rounded \\
\hline 
comp01.course & 4.00 & 4 s & 5.00 & 0 s & 7 \\ 
comp02.course & 5.98 & 2 s & 6.00 & 7 s & 12 \\ 
comp03.course & 6.94 & 1 s & 7.00 & 7 s & 14 \\ 
comp04.course & 4.98 & 1 s & 5.00 & 3 s & 12 \\ 
comp05.course & 8.00 & 1 s & 7.99 & 3 s & 14 \\ 
comp06.course & 5.99 & 3 s & 6.00 & 8 s & 14 \\ 
comp07.course & 6.00 & 5 s & 6.55 & 21 s & 17 \\ 
comp08.course & 6.99 & 2 s & 6.98 & 7 s & 11 \\ 
comp09.course & 4.99 & 1 s & 5.00 & 5 s & 10 \\ 
comp10.course & 6.00 & 3 s & 6.39 & 11 s & 16 \\ 
comp11.course & 5.00 & 1 s & 6.00 & 0 s & 8 \\ 
comp12.course & 9.91 & 4 s & 9.96 & 15 s & 18 \\ 
comp13.course & 5.98 & 1 s & 6.00 & 7 s & 8 \\ 
comp14.course & 6.00 & 2 s & 6.00 & 10 s & 14 \\ 
comp15.course & 6.94 & 1 s & 7.00 & 7 s & 15 \\ 
comp16.course & 6.00 & 4 s & 5.99 & 9 s & 15 \\ 
comp17.course & 5.98 & 3 s & 6.00 & 14 s & 12 \\ 
comp18.course & 4.99 & 1 s & 5.22 & 1 s & 8 \\ 
comp19.course & 6.00 & 1 s & 6.00 & 5 s & 11 \\ 
comp20.course & 6.00 & 4 s & 6.37 & 14 s & 13 \\ 
comp21.course & 8.00 & 2 s & 8.00 & 14 s & 13 \\
\hline 
\end{tabular}

\end{table}

\subsection{Conflict Graphs from Timetabling Benchmarks}

As a further illustration of the strength of SDP lower bounds, we present lower bounds for conflict graphs from two timetabling benchmarks.
From instances used in Track 3 of International Timetabling Competition 2007, we have extracted course-based conflict graphs,
  where there is edge between two vertices, if there is a curriculum prescribing the enrollment in both corresponding courses,
  or if a single teacher teaches both courses. 
For details, please see \cite{Bonutti2010} or \cite{Marecek2008Patat}.
From Toronto Examination Timetabling Benchmark, we have extracted exam-based conflict graphs,
  where there is edge between two vertices, if there is a student who should sit both corresponding exams. 
For details, please see \cite{Qu2009}.
For each graph, we have computed the best possible vertex colouring, without any bound on the number of uses of a colour,
and took the size of the largest colour class to be $C$. 
Subsequently, we have obtained lower bounds, upper bounds, and optima for $(C-1)$-bounded colouring, $(C-2)$-bounded colouring, etc.
%The values are presented in Table~\ref{tab:Udine}.

%\subsection{Odd Holes}
%On the other hand, one can show there are graphs where theta-like bounds perform poorly. 
%The ``odd hole'' graphs are by the strong perfect graph theorem 
%the minimally imperfect graphs \cite{MR2233847}. For every odd $n \ge 5$, there is a single 
%``odd hole'' graph on vertices $1, 2, \ldots, n-1, n$.
%Its $n$ vertices are $\{ (1, 2), (2, 3), \ldots, (n-1, n), (n, 1) \}$. The largest clique is hence 2.

\begin{table}[t!]
\caption{For Knesser graphs $K(n, 2)$ and forbidden intersection graphs $F(n, \gamma)$,
 where the size of the largest colour class in an optimal colouring is $C$, 
lower bounds $\OUR^{m}$ and optima $\chi^{m}$ for $(C - m)$-bounded colouring are shown. For $m = 0$, no bounds were applied.}
\label{tab:Knesser}
%\begin{center}
\begin{tabular}{l|ll|ll|ll|ll}
Graph & $\OUR^{0}$ & $\chi^{0}$ & $\OUR^{-1}$ & $\chi^{-1}$ & $\OUR^{-2}$ & $\chi^{-2}$ & $\OUR^{-3}$ & $\chi^{-3}$ \\ 
\hline 
$K(5,2)$      & 2.50 & 3 & 3.33 & 4 & 5.00 & 5 & 10.00 & 1 \\ 
$K(6,2)$      & 3.00 & 4 & 3.75 & 4 & 5.00 & 5 & 7.50 & 8 \\ 
$K(7,2)$      & 3.50 & 5 & 4.20 & 5 & 5.25 & 6 & 7.00 & 7 \\ 
$K(8,2)$      & 4.67 & 6 & 5.60 & 6 & 7.00 & 7 & 9.33 & 10 \\ 
$FI(6, 0.50)$ & 2.00 & 2 & 2.03 & 3 & 2.06 & 3 & 2.13 & 3 \\ 
$FI(6, 0.67)$ & 6.40 & 7 & 7.11 & 8 & 8.00 & 8 & 9.14 & 10 \\ 
$FI(6, 0.83)$ & 2.00 & 2 & 2.03 & 3 & 2.06 & 3 & 2.13 & 3 \\ 
$FI(6, 1.00)$ & 2.00 & 2 & 2.00 & 2 & 2.06 & 3 & 2.13 & 3 \\ 
\hline 
\end{tabular}

\end{table}

\subsection{Two Examples of Theoretical Interest}

To illustrate the weakness of the bound on certain graphs, we present the results for Knesser graphs of Lov{\'a}sz \cite{MR514625} 
 and the ``forbidden intersections'' graphs of Frankl and R{\"o}dl \cite{MR871675}. 
Knesser graph $K(n,k)$, $n > k > 1$, has  
$\binom{n}{k}$ vertices, corresponding to subsets of $\{1, 2, \ldots ,n\}$ of cardinality $k$.
Two vertices are adjacent if the corresponding subsets are disjoint.
Lov{\'a}sz has shown \cite{MR514625} the chromatic number of $K(n,k)$ is exactly $n - 2k + 2$,
despite the fact $K(n, k)$ has no triangle for $n > 3k$. 
Similarly, forbidden intersections graph $F(m, \gamma)$, $m \ge 1, 0 < \gamma < 1$, such that $(1-\gamma)m$ is an even integer,   
  has $2^m$ vertices, corresponding
  to sequences of $m$ bits (zeros and ones). Two vertices are adjacent, if the
  corresponding sequences differ in precisely $(1 - \gamma) m$ bits.
It is known the theta bound of {L}ov\'asz and related semidefinite programming 
 relaxations of graph colouring perform poorly on both ``forbidden intersections'' \cite{545462} 
 and Knesser graphs \cite{Karger1998}: the lower bound
 is $O(1)$ as $n$ grows, whereas the actual chromatic number grow $O(n)$ with $n$.  
Table~\ref{tab:Knesser} shows the lower bound gets tighter as the bound on the number of 
 uses of a colour gets tighter.

It should be noted that there is a large difference between clique and chromatic numbers in both Knesser and forbidden intersection graphs,
  which makes them quite unlike conflict graphs encountered in timetabling applications.
Although semidefinite programming lower bounds for graph colouring are weak on these graphs, they do tighten,
  as the bound on the number of uses of colours tightens.    
Nevertheless, the proposed lower bound is far from tight, in the worst case.

\section{Conclusions}
\label{sec:conclusions}

This paper has explored the limits of representability of extensions of graph colouring in semidefinite programming (SDP).
SDP clearly provides some of the strongest known relaxations in timetabling.
In particular, relaxations of simple timetabling problems related to Lov{\'a}sz theta provide useful lower bounds on the number of periods 
 required in the timetable, considering the conflict graph, the number of rooms,
 the capacities of rooms and special equipment available therein,
 and a pre-assignment of certain events to certain periods.

In such low-dimensional SDP relaxations, the colour assignment is not represented directly, 
but only in terms of the classes of equivalence of nodes assigned the same colour. 
This is sufficient for simple timetabling problems as described above, and makes the representation naturally invariant under the permutation of the colours. 

In contrast, many objectives in timetabling refer to time-based patterns of activities, e.g., whether events should be on the same day or not. 
These are not invariant under ``colour permutations'' and so the ``same colour'' representation is no longer sufficient. 
The matrix variable will need to capture the assignment of events to rooms as well as periods, and hence be constrained 
so that there is only a single event in each room-period pair. This gives a constraint on the rank of the matrix variable, 
which can be relaxed in a SDP. Despite the higher dimension of such relaxations, relaxations of rank-minimisation have proven very successful in many other fields \cite{Fazel2004},
and may turn out to be applicable also in timetabling. 
Modelling further and progressively more complex problems in
  semidefinite programming, with particular focus on relaxations one can solve fast,
  offers ample space for future work. 

In theory, one may wonder whether the relaxations as the best one can obtain in polynomial time assuming the unique games conjecture \cite{Khot2005}.
One could also seek approximation results for the problems we describe, either for the relaxations and rounding procedures of  this paper, or for novel ones. For example, one could obtain so-called lifted relaxations, e.g., using the method of moments  of \cite{lasserre2015introduction} applied to the copositive formulation, and to analyse the rounding therein, as \cite{Bansal2016} have done for job-shop scheduling. 
These would be an important advances in our understanding of scheduling and timetabling. 

\bibliography{sdp}  

%include{appendix} 
  
%\paragraph{Acknowledgements}  
%This paper builds upon two brief abstracts by the present authors \cite{Marecek2010-patat,Marecek2012-patat}.
%Jakub's work has been funded in part by the European Union Horizon 2020 Programme (Horizon2020/2014-2020), under grant agreement no.\  688380.

\end{document}